\documentclass{article}[12pt]
\usepackage{epsfig,makeidx,amsmath,amsfonts,latexsym,subfigure}

\newtheorem{prop}{Proposition}[section]
\newtheorem{lemma}{Lemma}[section]

\begin{document}

\vspace*{2cm}

\begin{center}
\large\textbf{{Notes on}}
\end{center}

\begin{center}
\Huge\textbf{{RISK THEORY}}
\end{center}

\vspace{0.5cm}

\begin{center}
\textbf{\large{by}}
\end{center}

\begin{center}
\textbf{\large{Anders Martin-L\"{o}f, Anders Sk\"{o}llermo}}
\end{center}

\vspace{5cm}


\vspace{0.3cm}

\begin{center}
\textbf{\large{Mathematical Statistics}}\linebreak
\textbf{\large{Department of Mathematics}}\linebreak
\textbf{\large{Stockholm University}}\linebreak
\textbf{\large{2011}}
\end{center}

\thispagestyle{empty}

\vfill\eject

\mbox

\thispagestyle{empty}

\vfill\eject

\tableofcontents

\thispagestyle{empty}

\vfill\eject

\mbox

\thispagestyle{empty}

\vfill\eject

\setcounter{page}{1}

\section{Introduction}

\noindent Risk theory is the part of insurance mathematics that is
concerned with stochastic models for the flow of payments in an
insurance business. The purpose of an insurance is in general to
level out fluctuations in the cost for the policyholder and to
replace the often strongly varying cost with a more predictable
flow of payments. To achieve this, a large group of risks -- a
``collective" -- is created in which the costs of an individual
member can be highly stochastic, but where the total cost is
levelled out as a consequence of the law of large numbers.\medskip

\noindent In these lecture notes we will describe some basic
natural models for ``risk processes" and derive various types of
asymptotic laws for the fluctuations in the amount of loss. We
will also investigate how the fluctuations depend on variables
such as reserve capital, premium amount, reinsurance arrangements,
size of the collective and distribution of the included variables.
Models for both life and property insurance will be
considered.\medskip

\noindent One can distinguish between two different types of
risks: the insurance risk and the uncertainty concerning the
future returns from the collected reserve capital. These notes
will mainly be concerned with the former type of risk, which is
generally better known from a statistical point of view because it
changes slower over time so that observed losses can be expected
to be relevant in predicting future losses. Also, an important
difference between the risk types is that uncertainty in for
instance the development of the interest can not be levelled out
in the same way as the first type of risk, since it can not be
decomposed as a sum of many contributions, obeying the law of
large numbers. However it is of interest to model the influence of
both risk types and indeed the substantial development of finance
mathematics during the last years has resulted in several models
for financial risks. In recent research these models are combined
with models from traditional risk theory in an interesting way,
and new types of contracts are being analyzed.\medskip

\noindent Risk theory as a branch of probability has a long
tradition, particularly within Swedish insurance research. Some of
the models that we will be interested in were formulated already
in the beginning of the 20th century in works by Filip Lundberg
and Harald Cram\'{e}r, and the theory of ruin probabilities that
we will consider was developed in the 1930-50's by Cram\'{e}r,
Esscher, Segerdahl and Arfwedson among others. This research
inspired the development of the theory for stochastic processes,
and during the 1960-80's it has turned out that many problems in
queuing theory, storage theory and risk theory are closely related
and can be solved by the same methods. This has resulted in
several simplifications of the theory in that technically
complicated analytical methods have been replaced by probabilistic
techniques which are more intuitive. In these notes we will, as
far as possible, use these probabilistic methods.

\vfill\eject

\section{Stochastic models for the total amount of loss during a fixed period}

\noindent In risk theory there are two basic models for the amount
of loss in an insurance collective: the individual model and the
collective model. Both these models are described in this section.
We also derive approximations for tail probabilities for the
distribution of the total amount of loss.

\subsection{The individual risk model}

\noindent In this model we consider a (large) number of individual
policies - for instance we can think of whole life assurances -
that are in effect during, let's say, one financial year. For each
of the policies there is a (small) probability $p_i$ that a loss
occurs, and a probability $q_i=1-p_i$ that no loss occurs. If a
loss occurs the amount $x_i$ is payed to the policyholder, where
$x_i$ is specified in the agreement. The losses are assumed to be
independent. Let $\{M_i\}$ be independent Bernoulli-variables with
$P(M_i=1)=1-P(M_i=0)=p_i$. Then the individual amount of loss can
be written as $x_iM_i$ and the total loss is given by
$X:=\sum_ix_iM_i$. Since the total loss is a sum of independent
random variables, it is natural to define its distribution via the
generating function E$[e^{\xi X}]$, which is the product of the
individual generating functions, that is,

\begin{eqnarray*}
\textrm{E}\left[e^{\xi X}\right] & = & \prod_i\textrm{E}\left[e^{\xi x_iM_i}\right]\\
 & = & \prod_i\left(q_i+p_ie^{\xi x_i}\right).
\end{eqnarray*}

\noindent The mean and variance of the individual losses are
E$[x_iM_i]=x_ip_i$ and Var$(x_iM_i)=x_i^2p_iq_i$, implying that
E$[X]=\sum_ix_ip_i$ and Var$(X)=\sum_ix_i^2p_iq_i$. Now, since $X$
is a sum of independent random variables, a natural approach might
be to approximate its distribution with a normal distribution with
these parameters, that is, one could believe that
$$
P\left(\frac{X-\textrm{E}[X]}{\sqrt{\textrm{Var}(X)}}>x\right)\approx
1-\Phi(x).
$$
However, this approximation often turns out to be quite poor
because of the fact that $p_i$ is typically very small so that
rather few losses occur even when the number of policies is large.
In such a situation it is more natural to approximate the
distribution of $X$ with a so called \emph{compound Poisson
distribution}, which is constructed as follows: Let $\{N_i\}$ be
independent Poisson distributed variables with E$[N_i]=\lambda_i$,
that is,
$$
P(N_i=n)=\frac{\lambda_i^n}{n!}e^{-\lambda_i}.
$$
Pick $\lambda_i$ so that $P(N_i=0)=q_i$ and put
$$
M_i=\left\{\begin{array}{ll}
                      0 & \mbox{if }N_i=0,\\
                      1 & \mbox{if }N_i\geq 1
                    \end{array}
            \right.
$$
Then $P(M_i=0)=1-P(M_i=1)=q_i$ and hence $M_i$ has the right
distribution. Moreover, when $p_i$ is small, $M_i=N_i$ with large
probability. To see this, note that

\begin{eqnarray*}
P(M_i\neq N_i) & = & P(N_i\geq 2)\\
 & = & 1-P(N_i=0)-P(N_i=1)\\
 & = & 1-e^{-\lambda_i}-\lambda_ie^{-\lambda_i}\\
 & \approx & 1-(1-\lambda_i+\lambda_i^2/2)-\lambda_i(1-\lambda_i)\\
 & = & \lambda_i^2/2.
\end{eqnarray*}

\noindent By the choice of $\lambda_i$, we have
$1-p_i=e^{-\lambda_i}$, and, since $e^{-\lambda_i}\approx
1-\lambda_i$, it follows that $p_i\approx \lambda_i$. Hence
$P(M_i\neq N_i)\approx p_i^2/2$ when $p_i$ is small. In this
situation it is natural to approximate $X$ with $S:=\sum_ix_iN_i$.
This quantity has a compound Poisson distribution and, since
$$
P(M_i\neq N_i\textrm{ for some }i)\leq \sum_iP(M_i\neq
N_i)=O\big(\sum_ip_i^2\big),
$$
the approximation is good if $\sum_ip_i^2$ is small.\medskip

\noindent Just as the distribution of $X$, the distribution of $S$
can be defined via its generating function. Remember that the
generating function for a Poisson distributed variable is given by

\begin{eqnarray*}
\textrm{E}[e^{\xi N_i}] & = & \sum_{n=0}^\infty e^{\xi
n}\frac{\lambda_i^n}{n!}e^{-\lambda_i}\\
& = & \exp\left\{\lambda_i\left(e^{\xi}-1\right)\right\}.
\end{eqnarray*}

\noindent Since $\{N_i\}$ are independent, we have

\begin{eqnarray*}
\textrm{E}[e^{\xi S}] & = & \prod_i\textrm{E}\left[e^{\xi
x_iN_i}\right]\\
& = & \prod_i\exp\left\{\lambda_i\left(e^{\xi x_i}-1\right)\right\}\\
& = & \exp\bigg\{\sum_i\lambda_i\left(e^{\xi x_i}-1\right)\bigg\}.
\end{eqnarray*}

\noindent Introduce the notation
$g(\xi)=\sum_i\lambda_i\left(e^{\xi x_i}-1\right)$. We then have
E$[e^{\xi S}]=e^{g(\xi)}$ or, equivalently, $g(\xi)=\log
\textrm{E}\left(e^{\xi S}\right)$. In what follows we will derive
approximations for the distribution of $S$ and thereby hopefully
also for the distribution of $X$. In this context it is worth
noting that $M_i\leq N_i$ for all $i$ and hence $X\leq S$ if
$x_i>0$ for all $i$. This implies that $P(X>x)\leq P(S>x)$ and
thus, if we can find an upper bound for $P(S>x)$, then this bound
is valid also for $P(X>x)$.\medskip

\begin{figure}
\begin{center}
\epsfig{file=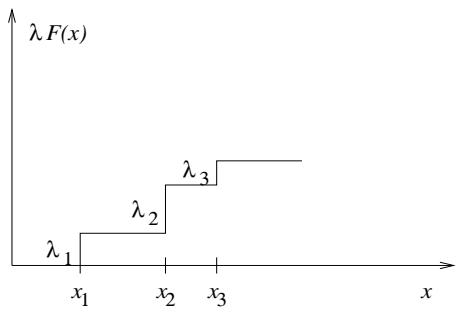,width=6cm,height=4cm} \caption{Construction
of $F(dx)$.} \label{fig:s8}
\end{center}
\end{figure}

\noindent The function $g(\xi)$ can be expressed in a slightly
different way using the so called \emph{risk mass distribution},
$F(dx)$. Let $\lambda=\sum_i\lambda_i$ and construct $F(dx)$ by
placing the mass $\lambda_i/\lambda$ at the point $x_i$ on the
x-axis, $i=1,2,3\ldots$, as demonstrated in Figure \ref{fig:s8}.
We then have

\begin{equation}\label{eq:Fog}
g(\xi)=\lambda\int_0^\infty\left(e^{\xi
x}-1\right)F(dx)\quad\textrm{and}\quad \textrm{E}\left[e^{\xi
S}\right]=e^{g(\xi)}.
\end{equation}

\noindent More generally we can consider $g(\xi)$ and $S$ defined
in this way with an arbitrary probability distribution $F(dx)$ and
some constant $\lambda<\infty$. The distribution of $S$ is then
called a compound Poisson distribution and the following
proposition gives a fundamental characterization for $S$.

\begin{prop}\label{prop:comppo}
Let $\{X_k\}$ be independent random variables with distribution
$F(dx)$ and let $N$ be a Poisson distributed variable, independent
of $\{X_k\}$, with E$[N]=\lambda$. Define $S=\sum_{k=1}^NX_k$.
Then $S$ has a compound Poisson distribution defined by
(\ref{eq:Fog}).
\end{prop}

\noindent \emph{Proof:} Let $f(\xi)$ be the generating function of
the distribution $F(dx)$, that is,
$$
f(\xi)=\textrm{E}\left[e^{\xi X_k}\right]=\int_0^\infty e^{\xi
x}F(dx).
$$
For each fixed $n$, the sum $S_n:=\sum_1^nX_k$ has distribution
$F^{n*}(dx)$ -- the convolution of $F$ with itself $n$ times --
with generating function E$\left[e^{\xi S_n}\right]=f^n(\xi)$.
Hence, $P(S\in dx|N=n)=F^{n*}(dx)$, and, summing over the possible
values of $N$, we obtain

\begin{equation}\label{eq:Sford}
P(S\in dx)=
\sum_{n=0}^\infty\frac{\lambda^n}{n!}e^{-\lambda}F^{n*}(dx).
\end{equation}

The corresponding generating function is

\begin{eqnarray*}
\textrm{E}\left[e^{\xi S}\right] & = & \sum_{n=0}^\infty
\textrm{E}\left[e^{\xi S}|N=n\right]P(N=n)\\
& = & \sum_{n=0}^\infty f^n(\xi)\frac{\lambda^n}{n!}e^{-\lambda}\\
& = & e^{\lambda(f(\xi)-1)}.
\end{eqnarray*}

\noindent With $g(\xi)$ defined as in (\ref{eq:Fog}), we have
$\lambda(f(\xi)-1)=g(\xi)$ and hence E$\left[e^{\xi
S}\right]=e^{g(\xi)}$, as desired. \hfill$\Box$

\subsection{The collective risk model}

\noindent In the individual risk model for a portfolio of whole
life assurances, the collective is changed over time as more and
more policyholders die. However, for moderate times and large
collectives this effect can often be neglected. A natural
approximation then is to consider a collective that is stationary
in time in the sense that $\lambda$ and $F(dx)$ are constant and
the number of losses in a time interval of length $t$ is Poisson
distributed with expected value $\lambda t$, the number of losses
in disjoint time intervals being independent. Below we give a
description of the total loss process $S(t)$ in the interval
$(0,t]$ motivated by this observation.\medskip

\begin{figure}
\begin{center}
\epsfig{file=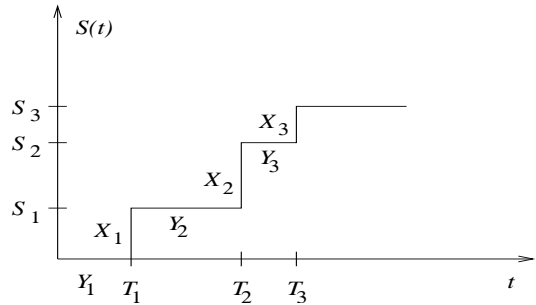,width=7cm,height=4cm} \caption{The total
loss process $S(t)$.}\label{fig:s12}
\end{center}
\end{figure}

\noindent Assume that the losses occur at time points
$T_1,T_2,\ldots$ that constitute a Poisson process in time, that
is, the increments $Y_k:=T_k-T_{k-1}$ are independent and
exponentially distributed with density $\lambda e^{-\lambda y}dy$.
At each time of loss $T_k$, an amount of damage $X_k>0$ is
generated. The variables $\{X_k\}$ are assumed to be independent
with distribution $F(dx)$ and the total loss in $(0,t]$ is given
by $S(t):=\sum_{T_k\in(0,t]}X_k$. As illustrated in Figure
\ref{fig:s12}, the process $S(t)$ is a step function with jumps of
height $X_k$ at the times $T_k$.\medskip

\noindent To specify the distribution of $S(t)$, let $N(t)$ denote
the number of losses in the interval $(0,t]$. We then have
$S(t)=\sum_{k=1}^{N(t)}X_k$. The process $\{N(t)\}_{t>0}$ is a
Poisson process with independent increments in disjoint intervals
and hence the increments of $S(t)$ -- that is, the sums of the
amounts of loss in disjoint intervals -- are also independent.
Furthermore, since
$$
P(N(t)=n)=\frac{(\lambda t)^n}{n!}e^{-\lambda t},
$$
by proceeding as in the derivation of (\ref{eq:Sford}), we obtain
$$
P(S(t)\in dx)=\sum_{n=0}^\infty \frac{(\lambda
t)^n}{n!}e^{-\lambda t}F^{n*}(dx).
$$
This means that, just like $S$ in the previous subsection, $S(t)$
has a compound Poisson distribution and hence its generating
function is given by

\begin{eqnarray*}
\textrm{E}\left[e^{\xi S(t)}\right] & = & \sum_{n=0}^\infty
\frac{(\lambda t)^n}{n!}e^{-\lambda t}f^n(\xi)\\
& = & e^{\lambda t(f(\xi)-1)}\\
& = & e^{tg(\xi)},
\end{eqnarray*}

\noindent where, as before, $f(\xi)$ is the generating function of
the distribution $F(dx)$ and
$g(\xi)=\lambda\int_0^\infty\left(e^{\xi x}-1\right)F(dx)$.
\medskip

\noindent The above formulas define the collective risk model,
which will be thoroughly studied in the following. The model can
be used to describe both a life assurance business and a property
insurance business. The total loss process $\{S(t)\}$ has
independent stationary increments with a compound Poisson
distribution defined by $\lambda$ and $F(dx)$ and the expected
value and variance of $S(t)$ can be obtained by differentiating
the generating function. Introducing the notation
$\mu=\int_0^\infty xF(dx)$ and $\nu=\int_0^\infty x^2F(dx)$, we
get
$$
\textrm{E}\left[S(t)\right]=tg'(0)=t\lambda\int_0^\infty
xF(dx)=t\lambda\mu
$$
and
$$
\textrm{Var}(S(t))=tg''(0)=t\lambda\int_0^\infty
x^2F(dx)=t\lambda\nu.
$$

\subsection{A method for calculating the distribution of $S(t)$}

\noindent Suppose that we have a fixed planning period. It is then
important to be able to calculate $P(S(t)>x)$ -- the probability
that the total loss exceeds $x$ -- as a function of $x$. In
general it is not possible to find simple formulas for this
probability. However, if the amounts of damage $X_k$ are
integer-valued -- that is, if $X_k\in\{1,2,3,\ldots\}$ -- then the
same thing holds for $S(t)$ and it turns out that we in this case
can derive a recursion formula for the masses of its distribution.
This so called \emph{Panjer-recursion} is easy to implement
numerically and is widely used. To describe it, assume for
simplicity that $t=1$ and write $S(1)=S$. Also, let
$f_x:=P(X_k=x)$ $(x=1,2,3,\ldots)$ and $g_y:=P(S=y)$
$(y=0,1,2,\ldots)$. Here the probabilities $\{f_x\}$ are assumed
to be known and we want to calculate $\{g_y\}$. To this end,
introduce the generating functions
$$
\varphi(s):=\sum_{x=1}^\infty s^xf_x\quad \textrm{and}\quad
\gamma(s):=\sum_{y=0}^\infty s^yg_y.
$$
Since $f(\xi)=\textrm{E}\left[e^{\xi X_k}\right]=\sum_x e^{\xi
x}f_x$, we have $f(\xi)=\varphi(e^\xi)$. Now let $\xi$ and $s$ be
related in that $s=e^\xi$. Then $f(\xi)=\varphi(s)$ and, since
$\gamma(s)=\textrm{E}\left[e^{\xi
S}\right]=e^{\lambda(f(\xi)-1)}$, we get
$$
\gamma(s)=e^{\lambda(\varphi(s)-1)}.
$$
Differentiating this relation we obtain
$\gamma'(s)=\lambda\varphi'(s)\gamma(s)$ or, more explicitly,

\begin{eqnarray*}
\gamma'(s) & = & \lambda\sum_{x=1}^\infty xf_xs^{x-1}
\sum_{y=0}^\infty g_ys^y\\
& = & \lambda \sum_{x=1}^\infty\sum_{y=0}^\infty
xf_xg_ys^{x+y-1}\\
& = & \lambda\sum_{n=1}^\infty s^{n-1}\sum_{x=1}^n xf_xg_{n-x}.
\end{eqnarray*}

\noindent But we also have $\gamma'(s)=\sum_{n=1}^\infty
ng_ns^{n-1}$. Equating these two expressions for $\gamma'(s)$
yields

\begin{equation}\label{eq:panje}
ng_n=\lambda\sum_{x=1}^nxf_xg_{n-x},\quad n=1,2,3,\ldots
\end{equation}

\noindent The probability $g_0$ is determined by noting that
$g_0=\gamma(0)=e^{\lambda(\varphi(0)-1)}=e^{-\lambda}$, where the
last equality follows since $\varphi(0)=0$. Given $g_0$, the
probabilities $\{g_n\}_{n\geq 1}$ are then successively obtained
from the equations (\ref{eq:panje}). We get

\begin{eqnarray*}
g_1 & = & \lambda f_1g_0 \\ g_2 & = & \lambda(f_1g_1+2f_2g_0)/2 \\
\vdots & &\vdots
\\g_n & = & \lambda(f_1g_{n-1}+2f_2g_{n-2}+\ldots+nf_ng_0)/n.
\end{eqnarray*}

\noindent As described above, an important quantity is
$G_m:=P(S>m)$, $m\geq 0$. Noting that $G_m=\sum_{m+1}^\infty g_n$,
the probabilities $\{G_n\}$ can be calculated together with
$\{g_n\}$ using the formula $G_m=G_{m-1}-g_m$, with $G_{-1}=1$.
Finally we remark that, if the $X_k$:s are not integer-valued,
they can be approximated by some suitable discretization and the
Panjer-recursion can then be applied to this distribution.

\subsection{Approximations of $P(S(t)>tx)$}

\noindent In this section we derive two useful approximations of
$P(S(t)>tx)$. They both involve the generating function $g(\xi)$
and are fairly easy to calculate when this function is
known.\medskip

\subsubsection{Chernoff bound}

\noindent The first approximation is based on an inequality,
\emph{Chernoff's inequality}, that is used in many statistical
contexts. To derive it, introduce the notation $F(t,dx)=P(S(t)\in
dx)$, fix $\xi\geq 0$, and note that

\begin{eqnarray*}
e^{tg(\xi)} & = & \textrm{E}\left[e^{\xi S(t)}\right]\\
& = & \int_0^\infty e^{\xi x}F(t,dx)\\
& \geq & e^{\xi t x}\int_{tx}^\infty F(t,dy)\\
& = & e^{\xi tx}P(S(t)\geq tx).
\end{eqnarray*}

\begin{figure}
\begin{center}
\subfigure[$g(\xi)$]{\epsfig{file=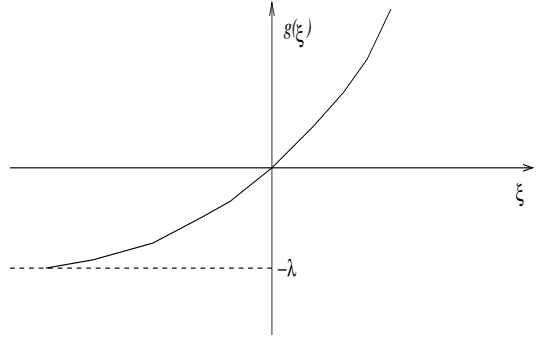,width=7cm,
height=4.5cm}}\par \subfigure[Construction of $h(x)$.]
{\epsfig{file=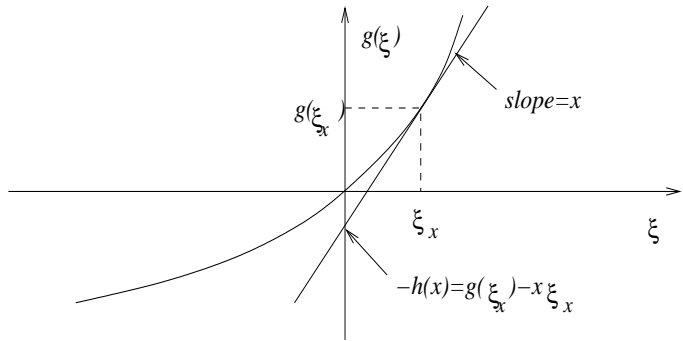,width=9cm,height=4.5cm}}\par
\subfigure[$h(x)$]{\epsfig{file=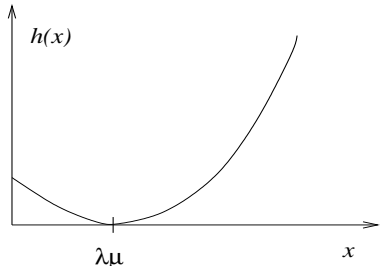,height=3.5cm,width=5cm}}
\caption{The functions $g(\xi)$ and $h(x)$.}\label{fig:s171820}
\end{center}
\end{figure}

\noindent Consequently we have $P(S(t)\geq tx)\leq
e^{-t(x\xi-g(\xi))}$ for all $\xi\geq 0$. Clearly the best upper
bound is obtained if $\xi\geq 0$ is picked so that $x\xi-g(\xi)$
is maximized. Define
$$
h(x)=\max_{\xi}\{x\xi-g(\xi)\}
$$
and write $\xi_x$ for the maximizing $\xi$-value. We then have
$$
P(S(t)\geq tx)\leq e^{-th(x)}\quad \textrm{if }\xi_x\geq 0.
$$
Analogously, it can be seen that
$$
P(S(t)\leq tx)\leq e^{-th(x)}\quad \textrm{if }\xi_x\leq 0.
$$
The function $h(x)$ will play an important role in what follows,
and we need to study its properties a bit closer. To this end,
first consider the function $g(\xi)=\lambda\int_0^\infty(e^{\xi
x}-1)F(dx)$. We will assume that $g(\xi)<\infty$ for
$\xi<\bar{\xi}$, where $\bar{\xi}>0$, that $g(\xi)\rightarrow
\infty$ as $\xi\rightarrow \bar{\xi}$ and also that
$g'(\xi)\rightarrow \infty$ as $\xi\rightarrow \bar{\xi}$. Since
$g'(\xi)=\lambda\int_0^\infty xe^{\xi x}F(dx)$ and
$g''(\xi)=\lambda\int_0^\infty x^2e^{\xi x}F(dx)$ are both
positive, the derivative $g'(\xi)$ increases monotonically from 0
to $\infty$ as $\xi$ increases from $-\infty$ to $\bar{\xi}$.
Hence $g(\xi)$ is strictly convex and increases from $-\lambda$ to
$\infty$ for these $\xi$-values, see Figure
\ref{fig:s171820}(a).\medskip

\noindent Now consider the function $h(x)$. The maximizing value
$\xi_x$ must satisfy $g'(\xi_x)=x$ and, since $g'(x)$ is strictly
increasing and continuous, for each $x$ this equation has exactly
one solution $\xi_x\leq \bar{\xi}$. Furthermore, the fact that
$g'(\xi)$ is strictly increasing also implies that $\xi_x\geq 0$
if and only if $x=g'(\xi_x)\geq g'(0)=\lambda \mu$. Hence
Chernoff's inequalities tells us that
$$
\textrm{(i)\hspace{0.2cm}}P(S(t)\geq tx)\leq e^{-th(x)}\quad
\textrm{if }x\geq \lambda\mu;
$$
$$
\textrm{(ii)\hspace{0.2cm}}P(S(t)\leq tx)\leq e^{-th(x)}\quad
\textrm{if }x\leq \lambda\mu.
$$
A picture of the geometrical construction of the function $h(x)$
is shown in Figure \ref{fig:s171820}(b). Consider the problem of
finding a tangent $-h+x\xi$, with given slope $x>0$, to the curve
$g(\xi)$. The tangent point $\xi_x$ satisfies $g'(\xi_x)=x$ and
$h=h(x)$ is determined so that $-h+x\xi_x=g(\xi_x)$, that is, we
have $h(x)=x\xi_x-g(\xi_x)$. As can be seen in the figure,
$h(x)\geq 0$ for all $x>0$ and $h(x)=0$ when $\xi_x=0$. The
geometrical construction can be thought of as if a line $-h+x\xi$,
with $x$ fixed, is pushed upwards towards the curve $g(\xi)$ until
a point is found where the line coincides with the tangent of the
curve. This means that we are looking for the smallest value of
$h$ such that $-h+x\xi\leq g(\xi)$ for all $\xi$, that is, such
that $h\geq x\xi-g(\xi)$ for all $\xi$. Hence the critical value
is $h(x)=\max_\xi\{x\xi-g(\xi)\}$.\medskip

\noindent The derivative of $h(x)$ is

\begin{eqnarray*}
h'(x) & = & \frac{d}{dx}(x\xi_x-g(\xi_x))\\
& = & \xi_x+\frac{d\xi_x}{dx}(x-g'(\xi_x))\\
& = & \xi_x,
\end{eqnarray*}

\noindent where the last equality follows because $x=g'(\xi_x)$.
The relation $x=g'(\xi_x)$ between $x$ and $\xi_x$ is 1-1 and
differentiable. We have $\frac{dx}{d\xi}=g''(\xi_x)$, and, since
$g''(\xi)>0$, it follows that $\frac{d\xi_x}{dx}=1/g''(\xi_x)$.
Using this, we get
$$
h''(x)=\frac{d\xi_x}{dx}=\frac{1}{g''(\xi_x)}>0,
$$
which means that $h(x)$ is also strictly convex. Remembering that
$h(x)\geq 0$ for all $x$, and $h(\lambda\mu)=0$, we can draw $h$
as in Figure \ref{fig:s171820}(c).\medskip

\noindent The relation between $g(\xi)$ and $h(x)$ can be
inverted. For $\xi=\xi_x$, we have $g(\xi)=x\xi-h(x)$ and
$\xi=h'(x)$. This implies that $g(\xi)$ is given by the formula
$g(\xi)=\max_x\{x\xi-h(x)\}$, which is analogous to the formula
for $h(x)$. In the theory for convex functions this relation is
well-known and $g(\xi)$ and $h(x)$ are said to be each others
\emph{Legendre transforms}.\medskip

\noindent Since $h(x)>0$ for $x\neq \lambda\mu$, Chernoff's
inequalities tells us that, when $x>\lambda\mu$ is fixed, the
probability of the event $\{S(t)/t\geq x\}$ decays exponentially
as $t\rightarrow \infty$. Such exponential estimates are common in
the theory of \emph{large deviations}. Here ``deviations" refer to
deviations from the mean and ``large" refers to the fact that the
deviations are large compared to the deviations treated by the
central limit theorem, where $x=\lambda\mu+y/\sqrt{t}$ as
$t\rightarrow \infty$. The function $h(x)$ that measures the decay
of the deviation probability is a fundamental object. It is called
the \emph{entropy function} of the distribution $F(dx)$. Let us
give some examples of how it is calculated for different
distributions.

\begin{itemize}
\item[1.] The exponential distribution, $F(dx)=e^{-x}dx$: For
$\xi<1$, we have

\begin{eqnarray*}
g(\xi) & = & \lambda\int_0^\infty(e^{\xi x}-1)e^{-x}dx\\
& = & \lambda\left(\frac{1}{1-\xi}-1\right)\\
& = & \frac{\lambda\xi}{1-\xi}.
\end{eqnarray*}

\noindent This yields $g'(\xi)=\lambda/(1-\xi)^2$ and the equation
$x=g'(\xi)$ hence becomes $1-\xi=\sqrt{\lambda/x}$. Thus

\begin{eqnarray*}
h(x) & = & x\xi-g(\xi)\\
& = & x\left(1-\sqrt{\frac{\lambda}{x}}\right)-
\lambda\left(\sqrt{\frac{x}{\lambda}}-1\right)\\
& = & x-2\sqrt{\lambda x}+\lambda\\
& = & \lambda\left(\sqrt{\frac{x}{\lambda}}-1\right)^2.
\end{eqnarray*}

\item[2.] The one-point distribution $F(dx)=\delta(x-1)dx$ gives
$g(\xi)=\lambda(e^{\xi}-1)$ and $g'(\xi)=\lambda e^{\xi}$. Putting
$x=g'(\xi)$, we get $\xi=\log(x/\lambda)$ and hence

\begin{eqnarray*}
h(x) & = & x\log\left(\frac{x}{\lambda}\right)-
\lambda\left(\frac{x}{\lambda}-1\right)\\
& = & \lambda\left[\frac{x}{\lambda}\log
\left(\frac{x}{\lambda}\right)-\frac{x}{\lambda}+1\right].
\end{eqnarray*}

\item[3.] The Gamma distribution with $\mu=a$,
$F(dx)=\gamma_a(x)dx$, gives $g(\xi)=\lambda((1-\xi)^{-a}-1)$ and
$g'(\xi)=\lambda a(1-\xi)^{-(a+1)}$. The relation $x=g'(\xi)$
implies that
$$
(1-\xi)=\left(\frac{\lambda a}{x}\right)^{1/(1+a)}
$$
and hence

\begin{eqnarray*}
h(x) & = & x\left(1-\left(\frac{\lambda
a}{x}\right)^{1/(a+1)}\right)-\lambda\left(\left(\frac{x}{\lambda
a}\right)^{a/(a+1)}-1\right)\\
& = &
\lambda\left[\frac{x}{\lambda}-\left(\frac{x}{\lambda}\right)^{a/(a+1)}
\left(a^{1/(a+1)}+a^{-a/(a+1)}\right)+1\right]
\end{eqnarray*}
\end{itemize}

\subsubsection{Esscher's approximation}

\noindent The second approximation of $P(S(t)>tx)$ is the so
called \emph{Esscher-approximation}, which is an asymptotic
formula, valid as $t\rightarrow \infty$. It states that
$$
P(S(t)>tx)\approx \frac{C}{\sqrt{t}}e^{-th(x)}\quad \textrm{as }
t\rightarrow\infty
$$
in the sense that the quotient between the left hand side and the
right hand side tends to 1. Here $C>0$ is a constant, and the
correction factor $C/\sqrt{t}$ gives a more precise estimate of
the exponential decay derived in the previous section. We will see
that in many cases this formula gives a good approximation also
for moderate values of $t$ and that it is easy to calculate
numerically if the function $g(\xi)$ is available.\medskip

\noindent Since the process $\{S(t)\}$ has independent increments,
it obeys the central limit theorem, that is,
$(S(t)-t\lambda\mu)/\sqrt{t\lambda\nu}$ is approximately normally
distributed as $t\rightarrow\infty$. This means that, with
$x=\lambda\mu+y\sqrt{\lambda\nu/t}$, we have that
$$
P(S(t)>tx)\rightarrow 1-\Phi(y)\quad \textrm{as
}\rightarrow\infty,
$$
where $\Phi(y)$ denotes the standard normal distribution function.
The central limit theorem hence gives an approximation for
``normal" deviations -- that is, deviations of the form
$y\sqrt{\lambda\nu/t}$ -- from the mean $\lambda\mu$. However, if
we want to study ``large" deviations, with $x>\lambda\mu$ fixed as
$t\rightarrow \infty$, then this approximation is not sufficient.
Below we will see that this problem can be circumvented by
modifying the distribution $F(dx)$ -- and thereby also the
distribution of $S(t)$ -- so that it becomes centered at the value
$tx$ that we are interested in. The central limit theorem can then
be applied to the transformed distribution to get an approximation
that can be used also for the original distribution close to the
value $tx$.\medskip

\noindent The modification of the distribution $F(dx)$ that we
will use is called the \emph{Esscher-transform}. It is obtained by
introducing a distribution that is proportional to $e^{ax}$ with
respect to $F(dx)$, where $a$ is a parameter that can be chosen
freely. To be more precise, we embed $F(dx)$ in an exponential
family by defining
$$
F_a(dx)=\frac{e^{ax}}{f(a)}F(dx),
$$
where $f(a)=\int_0^\infty e^{ax}F(dx)$. For $a<\bar{\xi}$ we have
$f(a)<\infty$ and hence $F_a(dx)$ is a probability distribution.
Now let $F_a(dx)$ be the modified distribution of $\{X_k\}$. It
defines a different distribution of $S_n=\sum_1^nX_k$. Write
$P_a(\cdot)$ for the modified probabilities and E$_a[\cdot]$ for
the corresponding means. Furthermore, let
$f_a(\xi):=\textrm{E}_a\left[e^{\xi X_1}\right]$ denote the
generating function of $F_a(dx)$. We then have

\begin{eqnarray}
f_a(\xi) & = & \int_0^\infty e^{\xi x}F_a(dx)\nonumber\\     
& = & \int_0^\infty e^{\xi x}\frac{e^{ax}}{f(a)}F(dx)\nonumber\\
& = & \frac{f(\xi+a)}{f(a)}.\label{eq:f_a}
\end{eqnarray}

\noindent Since the $X_k$:s are independent also under the measure
$P_a$, the distribution of $S_n$ under this measure is given by
$F_a^{n*}(dx)$ -- the convolution of $F_a$ with itself $n$ times.
Hence $\textrm{E}_a\left[e^{\xi S_n}\right]=\int_0^\infty e^{\xi
x}F_a^{n*}(dx)$. But, using (\ref{eq:f_a}), we also have

\begin{eqnarray*}
\textrm{E}_a\left[e^{\xi S_n}\right] & = & f_a^n(\xi)\\
& = & \frac{f^n(\xi+a)}{f^n(a)}\\
& = & \frac{\textrm{E}\left[e^{(a+\xi)S_n}\right]}{f^n(a)}\\
& = & \int_0^\infty e^{\xi x}\frac{e^{ax}}{f^n(a)}F^{n*}(dx).
\end{eqnarray*}

\noindent Thus
$$
F_a^{n*}(dx)=\frac{e^{ax}}{f^n(a)}F^{n*}(dx),
$$
that is, $F_a^{n*}$ is the Esscher-transform of $F^{n*}$. This
means that the original distribution of $S_n$ can be expressed in
terms of the modified one via the relation
$F^{n*}(dx)=f^n(a)e^{-ax}F_a^{n*}(dx)$. Hence, if we can
approximate $F_a^{n*}(dx)$ for some choice of $a$, we can also
approximate $F^{n*}(dx)$ via this relation. The reason for picking
an exponential density for $\{X_k\}$ is that this is the only case
when the transformed distribution of $S_n$ is obtained by applying
the same transform to the original distribution of $S_n$.\medskip

\noindent Now let us make an analogous transformation of $S(t)$.
Write $P_a(S(t)\in dx)=F_a(t,dx)$ and define
$$
F_a(t,dx)=e^{ax-tg(a)}F(t,dx).
$$
Remembering that E$\left[e^{\xi S(t)}\right]=e^{tg(\xi)}$, we then
have

\begin{eqnarray*}
\int_0^\infty F_a(t,dx) & = & e^{-tg(a)}\int_0^\infty e^{ax
}F(t,dx)\\
& = & e^{-tg(a)}\textrm{E}\big[e^{aS(t)}\big]\\
& = & 1
\end{eqnarray*}

\noindent so that $F_a(t,dx)$ is indeed a probability
distribution. The generating function is given by

\begin{eqnarray*}
\textrm{E}_a\big[e^{\xi S(t)}\big] & = & \int_0^\infty e^{\xi
x}e^{ax-tg(a)}F(t,dx)\\
& = & \textrm{E}\big[e^{(a+\xi)S(t)-tg(a)}\big]\\
& = & e^{t(g(a+\xi)-g(a))}.
\end{eqnarray*}

\noindent Hence we have E$_a\left[e^{\xi
S(t)}\right]=e^{tg_a(\xi)}$, where $g_a(\xi)=g(a+\xi)-g(a)$. This
means that $S(t)$ is still a compound Poisson process, since

\begin{eqnarray*}
g_a(\xi) & = & \lambda\int_0^\infty
\left(e^{(a+\xi)x}-e^{ax}\right)F(dx)\\
& = & \lambda f(a)\int_0^\infty \left(e^{\xi x}-1\right)F_a(dx).
\end{eqnarray*}

\noindent We thus have the important relation that, under the
measure $P_a$, $S(t)$ has a compound Poisson distribution with
$\lambda_a=\lambda f(a)$, jump distribution $F_a(dx)$ and
generating function $g_a(\xi)=g(a+\xi)-g(a)$. The last equation
immediately gives us the mean and variance. We have
$$
\textrm{E}_a[S(t)]=tg'_a(0)=tg'(a)
$$
and
$$
\textrm{Var}_a(S(t))=tg''_a(0)=tg''(a).
$$

\medskip

\noindent Just as for $S_n$, we have a simple expression for
$F(t,dx)$ in terms of $F_a(t,dx)$, namely
$$
F(t,dx)=e^{tg(a)-ax}F_a(t,dx).
$$
We will now see how this expression can be used to study large
deviations for $S(t)$. Consider the probability $P(S(t)\geq tx)$
with $x>\lambda\mu$. Center $P_a$ by choosing $a$ such that
E$_a[S(t)]=tx$, that is, such that $g'(a)=x$. We have previously
seen that this equation has a strictly positive unique solution if
$x>\lambda\mu=g'(0)$. The central limit theorem can now be used to
approximate the distribution of $S(t)$ under the measure $P_a$
near its mean $tx$: Put $S(t)=tg'(a)+Y$. Then, as $t\rightarrow
\infty$, the distribution of $Y$ is approximately normal with mean
0 and variance $\sigma^2=tg''(a)$. Furthermore,

\begin{eqnarray*}
P(S(t)\geq tx) & = & \int_{tx}^\infty F(t,dy)\\
& = & \int_{tx}^\infty e^{tg(a)-ay}F_a(t,dy)\\
& = & e^{tg(a)}\textrm{E}_a\left[e^{-a(tg'(a)+Y)},Y\geq 0\right]\\
& = & e^{t(g(a)-ag'(a))}\textrm{E}_a\left[e^{-aY},Y\geq 0\right].
\end{eqnarray*}

\noindent In the previous section we saw that, when $x=g'(a)$, we
have $ag'(a)-g(a)=h(x)$ and hence we arrive at the fundamental
formula
$$
P(S(t)\geq tx)=e^{-th(x)}\textrm{E}_a\left[e^{-aY},Y\geq 0\right].
$$
Now, if $Y$ has a density, the normal approximation for $Y$
implies that

\begin{eqnarray*}
\textrm{E}_a\left(e^{-aY},Y\geq 0\right) & \approx & \int_0^\infty
e^{-ay}\varphi\left(\frac{y}{\sigma}\right)\frac{dy}{\sigma}\\
& = & \int_0^\infty e^{-a\sigma y}\varphi(y)dy,
\end{eqnarray*}

\noindent where $\varphi(y)=e^{-y^2/2}/\sqrt{2\pi}$ denotes the
normal density. In the literature, 
\begin{eqnarray*}
E(s)&=& \int_0^\infty e^{-sy}\varphi(y)dy \\
       &=& e^{s^2/2}\int_0^\infty e^{-(s+y)^2/2}dy/\sqrt{2\pi} \\  
       &=& e^{s^2/2}(1-\Phi(s))    
 \end{eqnarray*} 
 is referred to as the \emph{Esscher function}
and, in terms of this function we have now derived Esscher's
approximation formula, which states that
$$
P(S(t)\geq tx)\approx e^{-th(x)}E(a\sigma)\quad \textrm{as
}t\rightarrow \infty,
$$
where $a>0$ is determined by the relation $x=g'(a)$ and
$\sigma=\sqrt{tg''(a)}$. The formula is only valid if $F(dx)$ has
a density, but later we will see that there is a similar
approximation if $F(dx)$ has a discrete distribution.\medskip

\noindent As $t\rightarrow\infty$, the same holds for $s$, and
from the definition of $E(s)$ we see that, for large $s$, the
exponential function is quickly damped as $y$ grows so that only
values near $y=0$ are essential. Near $y=0$, we have
$\varphi(y)\approx(1-y^2/2)/\sqrt{2\pi}$ and hence

\begin{eqnarray*}
E(s) & \approx & \frac{1}{\sqrt{2\pi}}\int_0^\infty
e^{-sy}\left(1-\frac{y^2}{2}\right)dy\\
& = & \frac{1}{\sqrt{2\pi}}\left(\frac{1}{s}-\frac{1}{s^3}\right).
\end{eqnarray*}

\noindent Thus we have the more explicit formula

\begin{equation}\label{eq:eappden}
P(S(t)\geq tx)\approx \frac{e^{-th(x)}}
{\sqrt{2\pi}a\sqrt{tg''(a)}},\quad x=g'(a),
\end{equation}

\noindent which is also referred to as Esscher's approximation.
The formula is reasonably easy to implement numerically provided
that it is possible to compute the function $g(a)$ and its
derivatives. If $x=g'(a)$, $h(x)=ag'(a)-g(a)$ and
$\sigma=\sqrt{tg''(a)}$ are computed for sufficiently many values
of $a>0$, the Esscher approximations can also be computed and
thereby we have an approximation for sufficiently many values of
$x$. This method gives an approximation that is good enough for
all distributions that occur in practice.\medskip

\noindent Even if the condition that $F(dx)$ has a density is not
fulfilled, it is possible to derive an analogous approximation
formula when $F(dx)$ is a discrete distribution such that $X_k$
takes values on the form $nd$, for some constant $d$ and
$n=0,1,2,\ldots$. To do this, note that if
$X_k\in\{nd;\hspace{0.1cm} n=0,1,2,\ldots\}$, the same thing holds
for $S(t)$. The normal approximation for $Y$ becomes
$$
P(Y=y)\approx \varphi\left(\frac{y}{\sigma}\right)
\frac{d}{\sigma}\quad\textrm{for }y=n\cdot d.
$$
and hence
$$
\textrm{E}\left[e^{-aY},Y\geq 0\right]\approx \sum_{n=0}^\infty
e^{-adn}\varphi\left(\frac{nd}{\sigma}\right)\frac{d}{\sigma}.
$$
In this case it is natural to introduce the discrete Esscher
function
$$
E(s,b)=\sum_{n=0}^\infty e^{-sn}\varphi(nb)b.
$$
The Esscher approximation then becomes
$$
P(S(t)\geq tx)\approx e^{-th(x)}E\left(ad,\frac{d}{\sigma}\right).
$$
As $b\rightarrow 0$ we have

\begin{eqnarray*}
E(s,b) & \approx & \sum_{n=0}^\infty e^{-sn}\varphi(0)b\\
& = & \frac{b}{\sqrt{2\pi}(1-e^{-s})},
\end{eqnarray*}

\noindent that is,
$$
E\left(ad,\frac{d}{\sigma}\right)\approx
\frac{1}{\sqrt{2\pi}\sigma}\left(\frac{d}{1-e^{-ad}}\right) \quad
\textrm{as }d\rightarrow \infty.
$$
Hence, in the discrete case we have the modified Esscher
approximation
$$
P(S(t)\geq tx)\approx
\frac{e^{-th(x)}}{\sqrt{2\pi}A(d)\sqrt{tg''(a)}},
$$
where $A(d)=(1-e^{-ad})/d$ and $x=g'(a)$. As $d\rightarrow 0$ we
see that $A(d)\rightarrow a$ and hence the formula is consistent
with (\ref{eq:eappden}).\medskip

\noindent The Esscher approximation holds analogously for
$P(S(t)\leq tx)$ when $x=g'(a)<\lambda\mu$ with $a<0$. More
generally, it holds for any probability $P(S(t)\in I)$, where $I$
is an interval $[z,y]$ with $\lambda\mu<z<y$ or $[y,z]$ with
$y<z<\lambda\mu$. In both cases, $a$ should be chosen so that
$g'(a)=x$, where $x$ is the point in $I$ where $h(x)$ is as small
as possible, that is, the exponent is always given by $\min_{x\in
I}h(x)$. The general formula is
$$
P(S(t)/t\in I)\approx\frac{C}{\sqrt{t}}e^{-t\min_{x\in I}h(x)}
$$
for some constant $C>0$. This type of estimate is common in the
more general theory for large deviations that has been developed
during the last decades inspired by the pioneering work of Esscher
from the 1930's.

\section{Theory of ruin probabilities}

\noindent So far we have studied the total loss $S(t)$ without
taking the flow of premiums in time into account, that is, we have
only considered $S(t)$ at a fixed time $t$. In this case it is
relevant to study $P(S(t)\ge tx)$ as we did in the previous
section. The number $tx$ should be thought of as the capital
available at time $t$ -- that is, the sum of the capital at $t=0$
and the amount of premiums that is paid in the interval $(0,t)$ --
and we want to make sure that this capital is large enough to make
the probability reasonably small. In such a setting we do not take
the possibility that a deficit might arise before time $t$ into
account.\medskip

\noindent To study the course of events in time we need to
describe the flow of premiums. This might also be stochastic, but
here we will restrict ourselves to the simplest setting, where the
premiums constitute a constant continuous inflow so that the total
premium paid in the interval $(0,t)$ is $ct$. If the capital at
time $t=0$ is $u$, the surplus at time $t$ is then given by
$u+ct-S(t)$ (for simplicity we disregard income from interest). In
the following we will study the so called \emph{ruin probability},
that is, the probability that the surplus is negative at some time
point during the planning period $(0,t)$, where $t=\infty$ is also
a possibility. In particular, we will see how this probability
depends on the parameters $u$, $c$, $\lambda$ and $F(dx)$.

\subsection{The total loss process}

\begin{figure}
\begin{center}
\epsfig{file=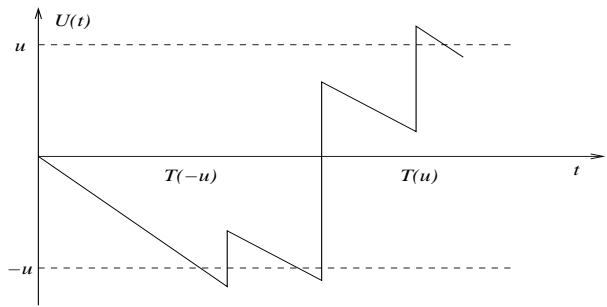,width=8cm,height=4cm} \caption{The net loss
cost process $U(t)$.}\label{fig:s36}
\end{center}
\end{figure}

\noindent Let us introduce the net amount of loss $U(t):=S(t)-ct$.
This is a stochastic process with upward jumps of height $\{X_k\}$
at times $\{T_k\}$, just as $S(t)$, and in between these times the
process decreases at rate $-c$; see Figure \ref{fig:s36}. Our main
object of interest is the time of ruin, denoted by $T(u)$ and
defined as the first time when $U(t)>u$. As we can see in Figure
\ref{fig:s36}, if $u\geq 0$, the ruin occurs at the first time
$T_k$ such that $S_k-cT_k>u$, that is, the ruin does not occur in
between two loss occasions, which means that in general a non-zero
deficit arises at time $T(u)$. We will also study $T(-u)$, which
is the first time when $U(t)\leq -u$. Since the heights of the
jumps are strictly positive, this occurs in between the jump
occasions so that, unlike what holds for $T(u)$, we have $U(t)=-u$
at time $t=T(-u)$. This will turn out to be a useful fact. Now
assume that $u\geq 0$ and define the ruin probabilities as

\begin{eqnarray*}
r(u,t) & = & P(T(u)\leq t),\quad \textrm{for }t<\infty;\\
r(u) & = & P(T(u)<\infty),
\end{eqnarray*}

\noindent and, analogously,

\begin{eqnarray*}
r(-u,t) & = & P(T(-u)\leq t),\quad \textrm{for }t<\infty;\\
r(-u) & = & P(T(-u)<\infty).
\end{eqnarray*}

\noindent We also define $T(\pm u)=\infty$ if the passage to $\pm
u$ never occurs, which, as we will see, happens with positive
probability.\medskip

\noindent Classical risk theory has to a large extent been
concerned with finding equations for $r(u)$ and $r(u,t)$ and, on
the basis of these equations, deriving approximations analogous to
the ones derived in the previous section for the distribution of
$S(t)$. In the following we will treat these problems, using more
probabilistic methods than the traditional ones. This often leads
to a better understanding of why the approximations are valid and
also to many simplifications of the derivations.\medskip

\noindent Before moving on to the mathematical treatment, we
remark that $T(u)$ is of course the natural ruin time when we have
a positive ``risk sum", that is, when the loss amounts $X_k$ are
positive and the premium inflow has rate $c>0$. This is the
natural model for property insurance and whole life assurance.
However, we can also apply the model to life assurance with
negative risk sum. In this case we have a continuous outflow of
payments $ct$ and $S(t)$ represents the accumulated inflow of
profits made at the times of the deaths. The ruin occurs when
$ct-S(t)\geq u$ for the first time, that is, at time $T(-u)$.
Hence $T(-u)$ also has a natural interpretation and $r(-u,t)$ and
$r(-u)$ are the ruin probabilities in this case.

\subsection{Basic formulas for the ruin probabilities}

\begin{figure}
\begin{center}
\epsfig{file=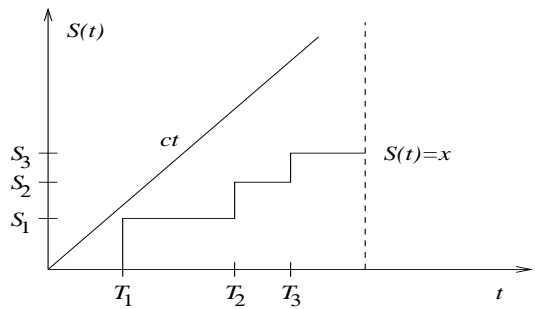,width=7cm,height=4cm} \caption{Illustration
of the event $A_t$.}\label{fig:s39}
\end{center}
\end{figure}

\noindent We begin by deriving a clever formula for the ruin
probability when $u=0$. The formula will turn out to be useful
also in finding expressions for the ruin probabilities when $u\neq
0$, as has been shown by Lajos Tak\'{a}cs. First consider the
event $A_t:=\{T(0)>t\}$ that the time to ruin exceeds $t$ and note
that

\begin{eqnarray*}
A_t & = & \{S(t')\leq ct'\textrm{ for all }t'\in(0,t)\}\\
& = & \{S_k\leq cT_k\textrm{ for }k=1,2,\ldots ,N(t)\},
\end{eqnarray*}

\noindent see Figure \ref{fig:s39} for an illustration. The
following lemma gives a simple formula for the probability of
$A_t$ given that $S(t)=x$, $0\leq x\leq ct$.

\begin{lemma}\label{lemma:ruin} We have
$$
P(A_t|S(t)=x)=\left(1-\frac{x}{ct}\right)_+,
$$
where
$$
\left(1-\frac{x}{ct}\right)_+=\left\{ \begin{array}{ll}
                      1-\frac{x}{ct} & \mbox{if\hspace{0.3cm}}0\leq x\leq ct,\\
                      0 & \mbox{if\hspace{0.3cm}}x>ct.
                    \end{array}
            \right.
$$
\end{lemma}

\noindent \emph{Proof:} We will use induction over $N(t)=n$ to
show the slightly stronger statement that

\begin{equation}\label{eq:rlemma}
P(A_t|S(t)=x,N(t)=n)=\left(1-\frac{x}{ct}\right)_+.
\end{equation}

\noindent To this end, first consider the case $n=0$. Then
$S(t)=0$, so that only $x=0$ has to be considered, and the event
$A_t$ occurs with probability 1. Hence (\ref{eq:rlemma}) is true
for $n=0$. For $n=1$, the event $A_t$ occurs if and only if
$T_1\geq x/c$. The conditional density for $T_1$ given that
$N(t)=1$ is

\begin{eqnarray*}
f_1(z)dz & = &
\frac{P(N(0,z)=0,N(z,z+dz)=1,N(z+dz,t)=0)}{P(N(0,t)=1)}\\
& = & \frac{e^{-\lambda z}\lambda dze^{-\lambda(t-z)}}{\lambda
te^{-\lambda t}}\\
& = & \frac{dz}{t},\qquad 0\leq z\leq t,
\end{eqnarray*}

\noindent that is, a uniform distribution on $(0,t)$. Hence
$$
P\left(T_1\geq \frac{x}{c}\hspace{0.1cm}\bigg|\hspace{0.1cm}
S(t)=x,N(t)=1\right)= \left(1-\frac{x}{ct}\right)_+,
$$
and so (\ref{eq:rlemma}) is true also for $n=1$.\medskip

\noindent Now assume that (\ref{eq:rlemma}) holds for $N(t)\leq
n-1$ and consider the case $S(t)=x$, $N(t)=n$. Given that
$N(t)=n$, the time $T_n$ has the conditional density $f_n(z)$
given by

\begin{eqnarray*}
f_n(z)dz & = & \frac{P(N(0,z)=n-1,N(z,z+dz)=1,N(z+dz,t)=0)}
{P(N(0,t)=n)}\\
& = & \frac{(\lambda z)^{n-1}e^{-\lambda z}\lambda
dze^{-\lambda(t-z)}/(n-1)!}{(\lambda t)^ne^{-\lambda t}/n!}\\
& = & n\left(\frac{z}{t}\right)^{n-1}\frac{dz}{t},\qquad 0\leq
z\leq t.
\end{eqnarray*}

\noindent If we fix $T_n=z$ and $S(T_{n-1})=y$, where $0\leq y\leq
x\leq cz\leq ct$, it follows from the induction assumption that
the conditional probability for $A_t$ is the same as for $A_z$,
that is, $1-y/cz$. Integrating over $y$ with the conditional
distribution of $S(T_{n-1})$ given $S(T_n)$ yields
$$
P(A_t|S(t)=x,N(t)=n,T_n=z)=\textrm{E}\left[1-\frac{S(T_{n-1})}{cz}\hspace{0.1cm}
\big|\hspace{0.1cm} S(T_n)=x\right].
$$
By symmetry we have
$$
\textrm{E}[S(T_{n-1})|S(T_n)]=(n-1)\textrm{E}[X_k|S(T_{n})]=\frac{n-1}{n}S(T_n),
$$
and hence
$$
P(A_t|S(t)=x,N(t)=n,T_n=z)=\left(1-\frac{(n-1)x}{ncz}\right).
$$
Integrating over $z\in(x/c,t)$ with the density $f_n(z)$ we
finally get

\begin{eqnarray*}
P(A_t|S(t)=x,N(t)=n) & = & \int_{x/c}^t\left(1-
\frac{(n-1)x}{ncz}\right)n\left(\frac{z}{t}\right)^{n-1}\frac{dz}{t}\\
& = & \int_{x/c}^tn\left(\frac{z}{t}\right)^{n-1}\frac{dz}{t}-
(n-1)\frac{xz^{n-2}}{ct^n}dz\\
& = & 1-\left(\frac{x}{ct}\right)^n-\frac{x}{ct}+
\left(\frac{x}{ct}\right)^n\\
& = & 1-\frac{x}{ct}.
\end{eqnarray*}

\noindent The formula (\ref{eq:rlemma}) now follows by induction.
Since the right hand side does not involve $n$, the conditioning
on $N(t)=n$ can be removed without affecting the formula and hence
the lemma is proved. \hfill$\Box$\medskip

\noindent Multiplying the probability in Lemma \ref{lemma:ruin}
with $P(S(t)\in dx)=F(t,dx)$ gives the joint probability
$$
P(A_t,S(t)\in dx)=\left(1-\frac{x}{ct}\right)_+F(t,dx).
$$
Now fix $S(t)=x$ such that $U(t)=S(t)-ct=x-ct\leq 0$ and write
$x-ct=-u$. We then have
$$
P(A_t|U(t)=-u)=\left(\frac{u}{ct}\right)_+
$$
and
$$
P(A_t,U(t)\in -du)=\left(\frac{u}{ct}\right)_+F(t,ct-du),
$$
that is,

\begin{equation}\label{eq:temp}
P(U(t)\in -du,U(t')\leq 0\textrm{ for }t'\in(0,t))=
\left(\frac{u}{ct}\right)_+F(t,ct-du);
\end{equation}

\noindent see Figure \ref{fig:s434445}(a). Integrating this over
$x\in (0,ct)$ we obtain the non-ruin probability
$\bar{r}(0,t):=1-r(0,t)$ for the initial capital $u=0$, that is,

\begin{equation}\label{eq:rbar}
\bar{r}(0,t)=\int_0^{ct}\left(1-\frac{x}{ct}\right)F(t,dx).
\end{equation}

\noindent In the following sections we will see how these formulas
can be used to determine the ruin probabilities when $u\neq 0$.

\subsubsection{The distribution of $T(-u)$}

\begin{figure}
\begin{center}
\subfigure[]{\epsfig{file=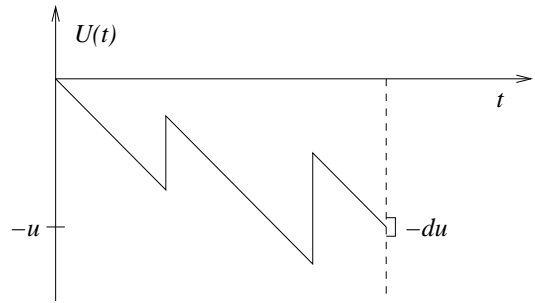,width=7cm, height=4cm}}\par
\subfigure[]{\epsfig{file=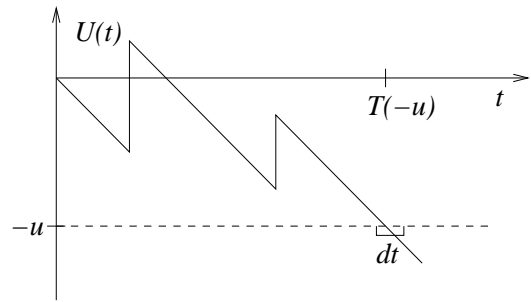,width=7cm,height=4cm}}\par
\subfigure[]{\epsfig{file=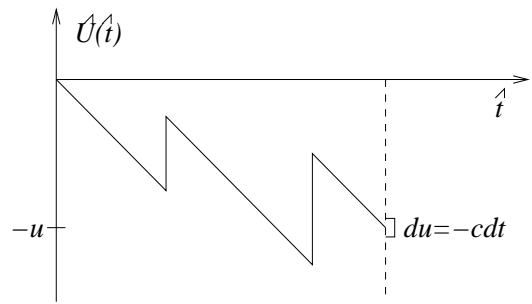,height=4cm,width=7cm}}
\caption{Transformation of $U(t)$.}\label{fig:s434445}
\end{center}
\end{figure}

\noindent Let us first derive a formula for $P(T(-u)\in dt)$. A
typical trajectory with $T(-u)\in dt$ has $U(t')>-u$ for
$t'<T(-u)$ and $U(t')=-u$ for $t'=T(-u)$; see Figure
\ref{fig:s434445}(b). If we turn this picture upside down and move
the origin to the crossing point, we see that the trajectory is
transformed into the trajectory in Figure \ref{fig:s434445}(a).
Hence, if this transformation does not change the distribution of
the process, the probability that $T(-u)\in dt$ should be the same
as the probability of the event in (\ref{eq:temp}), that is,
$$
P(T(-u)\in dt)=\left(\frac{u}{ct}\right)_+F(t,ct-du)
$$
and, since $-du=cdt$, we have
$$
P(T(-u)\in dt)=\left(\frac{u}{ct}\right)F(t,cdt-u)\quad
\textrm{for }ct\geq u>0.
$$
This is an explicit formula for the the distribution of $T(-u)$
and we have for instance that
$$
r(-u,t)=P(T(-u)\leq
t)=\int_{u/c}^t\left(\frac{u}{cs}\right)F(s,cds-u).
$$
To understand that the transformed process has the same
distribution as $U(t)$ we can write it as
$\hat{U}(\hat{t})=-u-U(t-\hat{t})$, $0\leq \hat{t}\leq t$. The
process $\hat{U}(\hat{t})$ has jumps at the time points
$\hat{T}_k=t-T_k$, which constitute a Poisson process, and the
jumps are $\hat{X}_k=X_k$, which are independent with distribution
$F(dx)$. Between the jumps, $\hat{U}(\hat{t})$ is changed at rate
$-c$. Hence $\hat{U}(\hat{t})$ is a process with the same
distribution as $U(t)$ and initial value $\hat{U}(0)=0$. The jumps
occur in a different order, but this does not affect the
distribution. The process $\hat{U}(\hat{t})$ is illustrated in
Figure \ref{fig:s434445}(c).\medskip

\noindent The ruin probability with $t=\infty$ is
$$
r(-u)=\int_{u/c}^\infty\left(\frac{u}{cs}\right)F(s,cds-u),
$$
or, with $x=cs-u$,
$$
r(-u)=\int_0^\infty\left(\frac{u}{x+u}\right)F\left(\frac{x+u}{c},dx\right).
$$
We will mainly consider the case when $c>\lambda\mu$ so that
E$[U(t)]=-(c-\lambda\mu)t<0$. By the law of large numbers,
$$
\frac{U(t)}{t}\rightarrow -(c-\lambda\mu)<0\quad\textrm{a.s. as
}t\rightarrow\infty.
$$
This implies that, with probability 1, the barrier $-u$ is hit
sooner or later, that is, $r(-u)=1$, or, equivalently,

\begin{equation}\label{eq:s46}
\int_0^\infty\left(\frac{u}{x+u}\right)F\left(\frac{x+u}{c},dx\right)=1\quad\textrm{if
}c>\lambda\mu.
\end{equation}

\noindent This relation will prove to be important in what
follows.

\subsubsection{The distribution of $T(u)$}

\noindent We will now derive an explicit formula for
$r(u,t)=P(T(u)\leq t)$ by using the previous results and
conditioning on the value of $U(t)$. Trivially
$$
r(u,t)=P(T(u)\leq t,U(t)>u)+P(T(u)\leq t,U(t)\leq u).
$$
If $U(t)>u$, we know for sure that $T(u)\leq t$, and hence

\begin{eqnarray*}
P(T(u)\leq t,U(t)>u) & = & P(U(t)>u)\\
 & = & \int_{x=u+ct}^\infty F(t,dx).
\end{eqnarray*}

\begin{figure}
\begin{center}
\epsfig{file=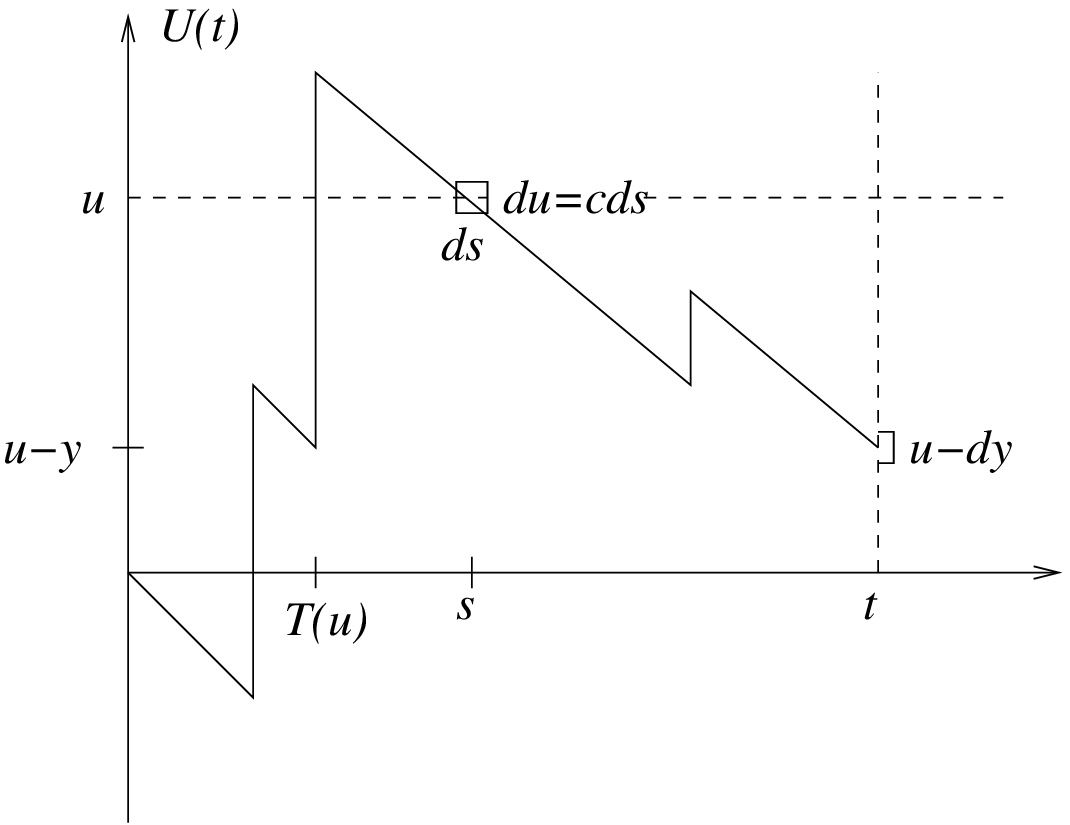,width=8cm,height=5cm} \caption{A scenario
with $U(t)\leq u$ and $T(u)\leq t$.}\label{fig:s47}
\end{center}
\end{figure}

\noindent If $T(u)\leq t$ and $U(t)\leq u$ the trajectory for
$U(t)$ has to cross the level $u$ one or more times between $T(u)$
and $t$; see Figure \ref{fig:s47}. Let $s$ be the value of the
last time when this occurs. The probability for such an outcome is
$P(U(s)\in du)P(E)$, where $E$ denotes the event to go from $u$ at
$s$ to $u-dy$ at $t$ without exceeding $u$ between $s$ and $t$.
The first factor equals $F(s,du+cs)=F(s,u+cds)$. By Lemma
\ref{lemma:ruin}, the last factor equals
$(y/c(t-s))F(t-s,c(t-s)-dy)$ and, integrating over $y\geq 0$, we
get $\bar{r}(0,t-s)$, see (\ref{eq:rbar}). Combining all this
yields
$$
r(u,t)=\int_{u+ct}^\infty F(t,dx)+\int_{0}^tF(s,u+cds)
\bar{r}(0,t-s).
$$
This formula is called \emph{Seals's formula} and, if $F(t,dx)$ is
known, it can be used to calculate $r(u,t)$. As $u$ and $t$
becomes large, it can also be used to derive an asymptotic formula
using the Esscher-approximation of $F(t,dx)$, but the calculations
become cumbersome.

\subsubsection{The ruin probability $r(u)$}

\begin{figure}
\begin{center}
\epsfig{file=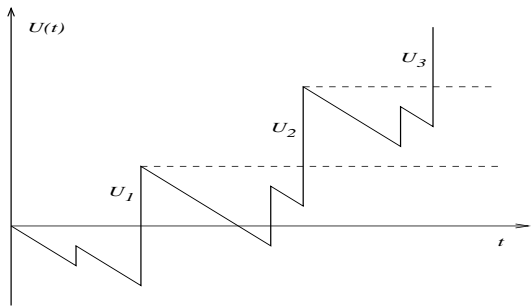,width=7cm,height=4cm} \caption{The process
$U(t)$ divided in up-crossings $\{U_k\}$.}\label{fig:s49}
\end{center}
\end{figure}

\noindent We will now derive a useful formula for
$r(u)=P(T(u)<\infty)$. A conceivable method for studying $r(u)$
would be to let $t\rightarrow\infty$ in Seal's formula for
$r(u,t)$. However, we will see that it is possible to obtain an
interesting formula via a more direct analysis, where the process
$U(t)$ is divided into successive \emph{upcrossings}, $\{U_k\}$;
see Figure \ref{fig:s49}. These upcrossings are defined as
follows: Initially, $U(0)=0$. With probability $r:=r(0)$, we have
$U(t)>0$ for some $t$, and with probability $1-r$, we have
$U(t)\leq 0$ for all $t$. In the first case, define $U_1$ to be
the value of $U(t)$ just after it has exceeded 0 for the first
time, that is, $U_1=U(T(0))$ if $T(0)<\infty$. From this point
$U(t)$ goes on for $t\geq T(0)$, and the process $U(T(0)+t)-U_1$
has the same distribution as $U(t)$ and is independent of $U_1$.
Define $U_2$ as the first up-crossing in this process, and so on.
In each step, there is a probability $1-r$ that no more
up-crossing occurs, and the successive $U_k$:s become independent
and identically distributed.\medskip

\begin{figure}
\begin{center}
\epsfig{file=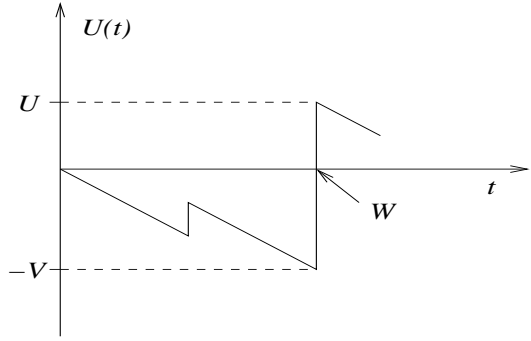,width=7cm,height=4.5cm} \caption{The
quantities $U$, $V$ and $W$.}\label{fig:s51}
\end{center}
\end{figure}

\noindent Below we will see that the $U_k$:s have a density $k(u)$
that is easy to write down and that
$$
r=\left\{\begin{array}{ll}
                      \frac{\lambda\mu}{c} & \mbox{if }c>\lambda\mu,\\
                      1 & \mbox{if }c<\lambda\mu.
                    \end{array}
            \right.
$$
Let $M$ denote the number of upcrossings and define
$\bar{U}=\sum_1^MU_k$. When $r<1$, $M$ has a geometric
distribution with
$$
P(M=m)=(1-r)r^m,\quad m=0,1,2,\ldots.
$$
This means that $M$ is finite with probability one and hence we
can write $\bar{U}=\max_{t\geq 0}U(t)$. The ruin probability then
becomes $r(u)=P(\bar{U}>u)$ and this probability can easily be
expressed in terms of $r$ and $k(u)$. It turns out, namely, that
$\bar{U}$ has a \emph{compound geometrical density},

\begin{equation}\label{eq:compge}
l(u)=(1-r)\sum_{m=0}^\infty r^mk^{m*}(u),
\end{equation}

\noindent where $k^{m*}(u)$ denotes the convolution of $k(u)$ with
itself $m$ times (compare with the compound Poisson distribution
characterized in Proposition \ref{prop:comppo}). This can be seen
by noting that, with probability $(1-r)r^m$, we have $M=m$ and the
density of $\bar{U}$ then becomes $k^{m*}(u)$. Summing over the
possible values of $m$, we get (\ref{eq:compge}). The formula for
$r(u)$ becomes
$$
r(u)=\int_u^\infty l(y)dy, \quad u>0.
$$
This formula is useful, since both $k(u)$ and $r$ can easily be
calculated. To find expressions for $k(u)$ and $r$, consider the
first upcrossing $U:=U_1$. Write $-V$ for the value of $U(t)$ just
before this up-crossing and let $W$ denote the time when the
up-crossing occurs; see Figure \ref{fig:s51}. By Lemma
\ref{lemma:ruin}, the joint distribution of $(U,V,W)$ is

\begin{eqnarray*}
P(U\in du,V\in dv, W\in dw) & = & P(U(w)\in -dv\hspace{0.1cm}
\textrm{and}\hspace{0.1cm}U(t)\leq 0\hspace{0.1cm}\textrm{for }t<
w)\cdot\\
& & P(\textrm{a loss occurs in }dw)\cdot\\
& & P(\textrm{the amount of loss}\in v+du)\\
& = & \left(\frac{v}{cw}\right)F(w,cw-dv)\lambda dwF(v+du).
\end{eqnarray*}

\noindent The distribution of $(U,V)$ is obtained by integrating
over $w$. Assuming that $F$ has a density $F'$ and making the
substitution $x=cw-v$, we get

\begin{eqnarray*}
P(U\in du, V\in dv) & = & \int_{w\geq
v/c}\left(\frac{v}{cw}\right) F'(w,cw-v)\lambda dw dv F'(u+v)du\\
& = & \lambda F'(u+v)dudv\int_0^\infty \left(\frac{v}{x+v}\right)
F'\left(\frac{x+v}{c},x\right)\frac{dx}{c}\\
& = & \frac{\lambda}{c}F'(u+v)dudv\quad \textrm{when }c\geq
\lambda\mu,
\end{eqnarray*}

\noindent where the last equality follows from (\ref{eq:s46}).
Hence the pair $(U,V)$ has density $(\lambda/c)F'(u+v)dudv$ for
$u,v\geq 0$. Integrating this over $v$ gives

\begin{eqnarray*}
P(U\in du) & = & \frac{\lambda}{c}du\int_0^\infty F'(u+v)dv\\
& = & \frac{\lambda}{c}(1-F(u))du.
\end{eqnarray*}

\noindent Since
$$
\int_0^\infty(1-F(u))du=\int_0^\infty uF'(u)du=\mu,
$$
we normalize by $\mu$, to get
$$
\frac{\lambda\mu}{c}\cdot\frac{1-F(u)}{\mu}du=rk(u)du,
$$
where $r=\lambda\mu/c<1$ and $k(u)=(1-F(u))/\mu$, with
$\int_0^\infty k(u)du=1$.\medskip

\noindent To summarize, we see that $P(U>0)=r=\lambda\mu/c$, and
that the conditional density of $U$ given that $U>0$, is given by
$k(u)=(1-F(u))/\mu$. Also, the ruin probability is
$$
r(u)=\int_u^\infty l(y)dy,\quad u>0,
$$
where $l(u)$ is specified in (\ref{eq:compge}). In risk theory,
this formula is called \emph{Cram\'{e}r's formula}, and in queuing
theory -- where it solves a similar problem -- it is referred to
as \emph{Pollaczek-Khinchin's formula}.\medskip

\noindent The formula for $l(u)$ gives rise to a corresponding
equation for the generating functions

\begin{equation}\label{eq:kappa}
\kappa(\xi)=\int_0^\infty e^{\xi u}k(u)du
\end{equation}

\noindent and

\begin{equation}\label{eq:lambda}
\lambda(\xi)=\int_0^\infty e^{\xi u}l(u)du,
\end{equation}

\noindent which are defined at least for $\xi\leq 0$. The equation
for $l(u)$ corresponds to the equation

\begin{eqnarray*}
\lambda(\xi) & = & (1-r)\sum_{m=0}^\infty r^m\kappa^m(\xi)\\
& = & \frac{1-r}{1-r\kappa(\xi)}.
\end{eqnarray*}

\noindent This equation can be used to compute $l(u)$ by inverting
the generating function.\medskip

\noindent \textbf{Example.} Let $k(u)=e^{-u}$. Then

\begin{eqnarray*}
\kappa(\xi) & = & \int_0^\infty e^{\xi u}e^{-u}du\\
& = & \frac{1}{1-\xi},
\end{eqnarray*}

\noindent which implies that

\begin{eqnarray*}
\lambda(\xi) & = & \frac{1-r}{1-r/(1-\xi)}\\
& = & (1-r)+\frac{r(1-r)}{1-r-\xi}.
\end{eqnarray*}

\noindent This is the generating function of
$l(u)=(1-r)\delta(u)+r(1-r)e^{-(1-r)u}$. The density $l(u)$ can be
computed in a similar way when $\kappa(\xi)$ is a rational
function of $\xi$.

\subsubsection{Panjer-approximation of $r(u)$}

\noindent In the following two sections, two different
approximations of $r(u)$ will be derived. The first one is a kind
of Panjer-recursion and the second one is an asymptotic formula as
$u\rightarrow\infty$ similar to the Esscher approximation.\medskip

\noindent As for the Panjer-recursion, consider first two discrete
distributions $\{k_n\}_1^\infty$ and $\{l_n\}_0^\infty$, where
$\{k_n\}$ is known and $\{l_n\}$ is ``compound geometric", that
is,
$$
l_n=(1-r)\sum_{m=0}^\infty r^mk^{m*}_n,\quad n=0,1,2,\ldots.
$$
We will now see that $\{l_n\}$ can be computed by aid of a
recursive formula of Panjer-type. Convolving the equation
$l=(1-r)\sum_0^\infty r^mk^{m*}$ with $rk$ yields

\begin{eqnarray*}
rk*l & = & (1-r)\sum_{m=0}^\infty r^{m+1}k^{(m+1)*}\\
 & = & (1-r)\sum_{m=1}^\infty r^mk^{m*}\\
 & = & l-(1-r)\delta,
\end{eqnarray*}

\noindent where
$$
\delta_n=k_n^{0*}=\left\{\begin{array}{ll}
                      1 & \mbox{if }n=0,\\
                      0 & \mbox{if }n>0.
                    \end{array}
            \right.
$$
Hence we have the renewal equation $l=(1-r)\delta+rk*l$, or, more
explicitly,
$$
l_n=(1-r)\delta_n+r\sum_{m=1}^nk_ml_{n-m}.
$$
The probabilities $\{l_n\}$ can successively be computed for
$n=0,1,2,\ldots$. We obtain

\begin{eqnarray*}
l_0 & = & 1-r\\
l_1 & = & rk_1l_0\\
l_2 & = & r(k_1l_1+k_2l_0)\\
\vdots & & \vdots\\
l_n & = & r(k_1l_{n-1}+\ldots +k_nl_0).
\end{eqnarray*}

\noindent If $\{k_n\}$ is known, this recursion is easy to
implement. Also, the probabilities $r_n=\sum_n^\infty l_m$ can be
computed parallel to $l_n$.\medskip

\noindent The equations for $l(u)$ and $r(u)$ are similar, but,
just as $k(u)$, they are the densities of continuous random
variables. However, if we make a suitable discrete approximation
of $k(u)$, we can calculate the corresponding approximations of
$l(u)$ and $r(u)$ by the above method.\medskip

\noindent The relation between $\{k_n\}$ and $\{l_n\}$ can also be
expressed in terms of the generating functions
$\hat{k}(s)=\sum_1^\infty k_ns^n$ and $\hat{l}(s)=\sum_0^\infty
l_ns^n$. We have

\begin{eqnarray*}
\hat{l}(s) & = & (1-r)\sum_{m=0}^\infty r^m\hat{k}^m(s)\\
 & = & \frac{1-r}{1-r\hat{k}(s)}.
\end{eqnarray*}

\noindent If, for instance, $\hat{k}(s)$ is a rational function of
$s$, the function $\hat{l}(s)$ is also rational and $\{l_n\}$ can
be obtained by partial fraction expansion.\medskip

\noindent \textbf{Example.} Let $k_n=(1-p)p^{n-1}$,
$n=1,2,\ldots$. This gives

\begin{eqnarray*}
\hat{k}(s) & = & (1-p)\sum_{n=1}^\infty s^np^{n-1}\\
 & = & \frac{(1-p)s}{1-ps},
\end{eqnarray*}

\noindent that is,

\begin{eqnarray*}
1-r\hat{k}(s) & = & 1-\frac{r(1-p)s}{1-ps}\\
& = & \frac{1-qs}{1-ps},
\end{eqnarray*}

\noindent where $q=p+r(1-p)$. Hence

\begin{eqnarray*}
\hat{l}(s) & = & (1-r)\frac{1-ps}{1-qs}\\
 & = & (1-r)\left(1+\frac{(q-p)s}{1-qs}\right)\\
 & = & (1-r)\bigg(1+r(1-p)\sum_1^\infty s^nq^{n-1}\bigg).
\end{eqnarray*}

\noindent This yields
$$
\left\{\begin{array}{l}
                      l_0=(1-r)\\
                      l_n=(1-r)r(1-p)q^{n-1},\hspace{0.1cm}n\geq 1.
                    \end{array}
\right.
$$
A natural method for finding a discrete approximation to the density $k(x)$ can be obtained as follows:
Approximate first the distribution $F(x)$ by a discrete distribution with masses $f_n$ for $x=nd$ and
$n=1,2, \ldots$ and put $F_n = \sum_1^n f_m$.
 For this distribution $F(x)$ is piecewise constant:
$F(x) = F_n$ and $1-F(x) =1-F_n = \sum_{n+1}^\infty f_m$ for $nd\leq x<(n+1)d,$ and then
$\mu = d\sum_0^\infty(1-F_n).$
The density $k(x) = (1-F(x))/\mu$ can then be approximated by a discrete distribution having masses
$k_n = \int_{(n-1)d}^{nd} k(x)dx = d(1-F_{n-1})/\mu = (d/\mu)\sum_n^\infty f_m$ for $x = nd$ and $n=1,\ldots.$. This distribution will have total mass one and is located at positive x-values.

\subsubsection{Cram\'{e}r-Lundberg's approximation of $r(u)$}

\noindent We will now derive a more explicit approximation formula
for $r(u)$. It is an asymptotic formula valid as $u\rightarrow
\infty$ and, as we will see, it is closely related to the Esscher
approximation.\medskip

\noindent First recall from Section 3.2.3 that
$$
r(u)=\int_u^\infty l(y)dy\quad \textrm{for }u>0,
$$
where $l(y)=(1-r)\sum_0^\infty r^mk^{m*}(y)$, $r=\lambda\mu/c$ and
$k(u)=(1-F(u))/\mu$. Here $r=P(U>0)$, where $U$ denotes the size
of an up-crossing, and $k(u)$ is the conditional density of $U$
given that $U>0$. To get an approximation of $l(y)$ when $y$ is
large, we need an approximation of $k^{m*}(y)$ as
$y\rightarrow\infty$. Since $r^m$ damps large $m$-values in the
formula for $l(y)$ we only have to consider moderate values of
$m$. The desired approximation is obtained by introducing a
modified density $k_a(y)$ as in Section 2.4.2, and choosing $a$
suitably. We have
$$
k_a(y)=\frac{e^{ay}k(y)}{\kappa(a)},
$$
where $\kappa$ is the generating function of the density $k$ (see
(\ref{eq:kappa})), and below we will see that, just as $g(a)$,
$\kappa(a)<\infty$ for $a<\bar{\xi}$. As in Section 2.4.2, we get
$$
k_a^{m*}=\frac{e^{ay}k^{m*}(y)}{\kappa^m(a)}
$$
so that
$$
k^{m*}(y)=e^{-ay}\kappa^m(a)k_a^{m*}(y)
$$
and hence
$$
l(y)=e^{-ay}(1-r)\sum_{m=0}^\infty r^m\kappa^m(a)k_a^{m*}(y).
$$
Choosing $a$ such that $r\kappa(a)=1$ yields
$$
l(y)=e^{-ay}(1-r)\sum_{m=0}^\infty k_a^{m*}(y).
$$
As $y\rightarrow \infty$, this expression can be approximated
using the so called \emph{renewal theorem}, which is an important
result in renewal theory. It states that, as $y\rightarrow\infty$,
the sum $\sum_0^\infty k_a^{m*}(y)$ can be approximated by a
uniform density with intensity $1/m_a$, where

\begin{eqnarray*}
m_a & = & \int_0^\infty yk_a(y)dy\\
 & = & \int_0^\infty \frac{ye^{ay}k(y)}{\kappa(a)}dy\\
 & = & \frac{\kappa'(a)}{\kappa(a)}.
\end{eqnarray*}

\noindent Substituting this approximation in the formula for
$r(u)$ gives

\begin{eqnarray*}
r(u) & \approx & \int_u^\infty e^{-ay}(1-r)\frac{dy}{m_a}\\
 & = & e^{-au}\frac{(1-r)}{m_a}\int_0^\infty e^{-ax}dx\\
 & = & \frac{(1-r)}{am_a}e^{-au}.
\end{eqnarray*}

\noindent If $a$ is known, this is a simple exponential
approximation.\medskip

\begin{figure}
\begin{center}
\epsfig{file=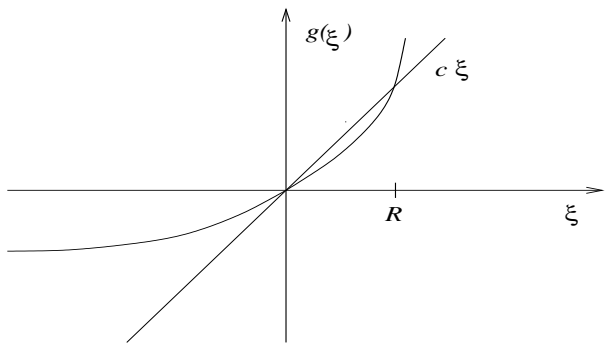,width=8cm,height=4.5cm}
\caption{Construction of $R$.}\label{fig:s60}
\end{center}
\end{figure}

\noindent The equation for $a$, $r\kappa(a)=1$, can be expressed
more explicitly in terms of $g(a)$, by noting that

\begin{eqnarray*}
\kappa(a) & = & \int_0^\infty
e^{ax}\left(\frac{1-F(x)}{\mu}\right)dx\\
& = & \frac{1}{a\mu}\int_0^\infty (e^{ax}-1)F(dx)\\
& = & \frac{g(a)}{a\lambda\mu},
\end{eqnarray*}

\noindent where the second equality is obtained by partial
integration. The equation for $a$ hence becomes
$g(a)/a\lambda\mu=1/r=c/\lambda\mu$, that is, $g(a)=ca$. Recall
from Section 2.4.1 that $g(\xi)$ is strictly convex with $g(0)=0$
and $g'(0)=\lambda\mu$. We are looking for the intersection with a
line $c\xi$ with slope $c$; see Figure \ref{fig:s60}. For
$c>g'(0)=\lambda\mu$, there is a strictly positive root, which is
denoted by $R$ and referred to as the \emph{Lundberg exponent}.
For $c<\lambda\mu$, the root is negative.\medskip

\noindent To find an expression for the constant $C:=(1-r)/am_a$
in the formula for $r(u)$, note that

\begin{eqnarray*}
m_a & = & \frac{\kappa'(a)}{\kappa(a)}\\
& = & [\log\kappa(a)]'\\
& = & [\log g(a)]'-[\log a]'\\
& = & \frac{g'(a)}{g(a)}-\frac{1}{a}\\
& = & \frac{g'(a)-c}{ca}.
\end{eqnarray*}

\noindent This yields

\begin{eqnarray*}
C & = & \frac{1-\lambda\mu/c}{(g'(a)-c)/c}\\
& = & \frac{c-g'(0)}{g'(a)-c}.
\end{eqnarray*}

\noindent To sum up, we have deduced that $r(u)\approx Ce^{-Ru}$,
where $R$ is the positive root of the equation $g(a)=ca$ and
$C=(c-g'(0))/(g'(R)-c)$. Here ``$\approx$" means that the quotient
between the right hand and the left hand side tends to 1 as
$u\rightarrow\infty$.
A natural way of using these formulas for the design of a system is to start by choosing c so that r has a suitable value close enough to one, and then finding the corresponding values of R and C. Then u can easily be found so that $Ce^{-Ru}$ has a value considered to be small enough to be safe.

\noindent \textbf{Example.} Approximate calculation of $R$ and $C$ when $r=\lambda\mu/c$ is close to one.
The equation $g(R)/R=c$ can be expressed in terms of the Taylor expansion of $g(R)$ as follows:
$$
g(R)=\lambda\sum_{k=1}^\infty\mu_kR^k/k!
$$
where $\mu_k$ is the k-th moment of the claims distribution F. ($\mu_1=\mu$ and $\mu_2=\nu$).
In terms of it the equation for R is hence
$$
\lambda(\mu+\mu_2R/2+\mu_3R^2/6+\cdots)=c
$$
or
$$
\lambda(\mu_2R/2+\mu_3R^2/6+\cdots)=\lambda\mu(1/r-1),
$$
and we see that $r\approx1$ corresponds to $\rho\equiv(1/r-1)\approx0$ and hence to $R\approx0$.
To first order in $\rho$ we hence have $\mu_2R_1/2=\mu_1\rho $ and $ R_1=(2\mu_1/\mu_2)\rho$.
To second order in $ \rho$ we then have $\mu_2R_2/2+\mu_3R_1^2/6=\mu_1\rho$ and
$$
R_2=R_1-(2/\mu_2)(\mu_3/6)(2\mu_1/\mu_2)^2\rho^2=(2\mu_1/\mu_2)\rho-(4/3)(\mu_3\mu_1^2/\mu_2^3)\rho^2
$$
etc.
The corresponding values of C can be obtained from the relation
\begin{eqnarray*}
C &=& (g(R)/R-\lambda\mu)/(g'(R)-g(R)/R)\\
    &=& (\sum_{k=2}^\infty\mu_kR^{k-1}/k!)/(\sum_{k=2}^\infty\mu_kR^{k-1}(k-1)/k!)\\
    &=& (\sum_{k=2}^\infty\mu_kR^{k-2}/k!)/(\sum_{k=2}^\infty\mu_kR^{k-2}(k-1)/k!)\\
    &=&\frac{(\mu_2/2)+(\mu_3/6)R+(\mu_4/24)R^2+\cdots}{(\mu_2/2)+(\mu_3/3)R+(\mu_4/8)R^2+\cdots}
 \end{eqnarray*}
 To first order in $\rho$ we hence have
 $
 C_1= ((\mu_2/2)+(\mu_3/6)R_1)/((\mu_2/2)+(\mu_3/3)R_1),
 $
 and we get quite explicit expressions in terms of the moments $\mu_k$ and $\rho$.

\begin{figure}
\begin{center}
\epsfig{file=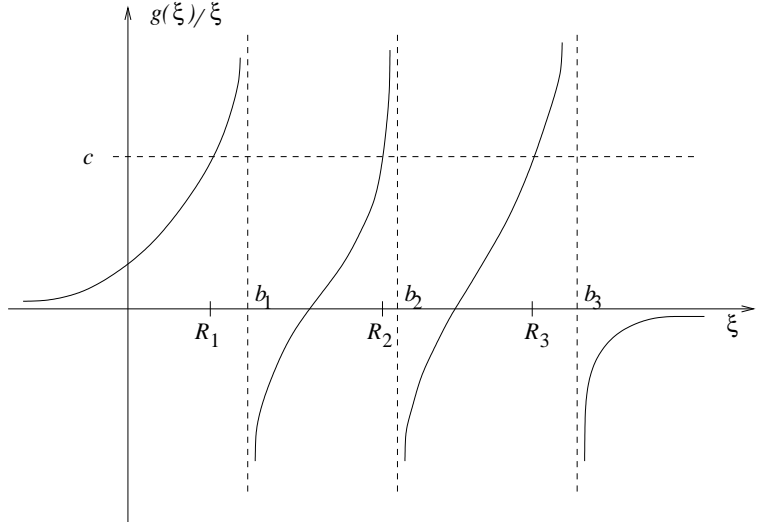,width=10cm,height=7cm} \caption{Graphical
picture of the roots $\{R_i\}$.}\label{fig:s61b}
\end{center}
\end{figure}

\noindent \textbf{Example.} Assume that $F'(x)$ is a weighted sum
of exponential densities, that is,
$$
F'(x)=\sum_{i=1}^na_ib_ie^{-b_ix},
$$
where $a_i>0$, $\sum_1^na_i=1$ and $0<b_1<b_2<\ldots<b_n$. Then
$r(u)$ can be calculated fairly explicitly via the generating
functions $\kappa(\xi)$ and $\lambda(\xi)$, defined in
(\ref{eq:kappa}) and (\ref{eq:lambda}) respectively, and we will
be able to see how the approximation $Ce^{-R\mu}$ arises. First
recall that $\kappa(\xi)=g(\xi)/\xi\lambda\mu$ and
$\lambda(\xi)=(1-r)/(1-r\kappa(\xi))$. The generating function,
$f(\xi)$, of $F'(x)$ is
$$
f(\xi)=\sum_{i=1}^na_i\cdot \frac{b_i}{b_i-\xi}
$$
and we obtain

\begin{eqnarray}
g(\xi) & = & \lambda(f(\xi)-1)\nonumber\\
& = & \lambda\sum_{i=1}^n\frac{a_i\xi}{b_i-\xi}\label{eq:exg}
\end{eqnarray}

\noindent and $\mu=\sum_1^na_i/b_i$. Hence
$$
\kappa(\xi)=\frac{1}{\mu}\sum_{i=1}^n\frac{a_i\xi}{b_i-\xi},
$$
that is, $\kappa(\xi)$ is a rational function of $\xi$, where the
denominator is of degree $n$, and $\kappa(\xi)\rightarrow 0$ as
$|\xi|\rightarrow \infty$. The poles of $\lambda(\xi)$ -- that is,
the zeroes of its denominator -- are the roots of the equation
$1-r\kappa(\xi)=0$. If the root $\xi=0$ is ignored, this equation
can be rewritten as $1=rg(\xi)/\lambda\mu \xi$, that is,
$g(\xi)=c\xi$. Using the relation (\ref{eq:exg}), the equation
becomes
$$
\lambda\sum_{i=1}^n\frac{a_i}{b_i-\xi}=c.
$$
For $\xi=0$, the left hand side equals $\lambda\mu$ which is
strictly smaller than $c$. A graph of the expression on the left
hand side as a function of $\xi\geq 0$ is displayed in Figure
\ref{fig:s61b}. We see that there are $n$ real roots $R_1,\ldots
R_n$, with $0<R_1<b_1<R_2<b_2<\ldots<R_n<b_n$. Hence, the partial
fraction expansion of $\lambda(\xi)$ is
$$
\lambda(\xi)=(1-r)+\sum_{i=1}^n\frac{C_iR_i}{R_i-\xi},
$$
where the coefficients $C_iR_i$ are determined by the formula

\begin{eqnarray*}
C_iR_i & = & \lim_{\xi\rightarrow R_i}(R_i-\xi)\lambda(\xi)\\
& = & \lim_{\xi\rightarrow
R_i}\frac{(R_i-\xi)(1-r)}{1-rg(\xi)/\xi\lambda\mu}\\
& = & \lim_{\xi\rightarrow
R_i}\frac{(R_i-\xi)(1-r)c\xi}{c\xi-g(\xi)}.
\end{eqnarray*}

\noindent Near $\xi=R_i$, we have for the denominator, that

\begin{eqnarray*}
g(\xi)-c\xi & \approx & g(R_i)-cR_i+(\xi-R_i)(g'(R_i)-c)\\
& = & (\xi-R_i)(g'(R_i)-c),
\end{eqnarray*}

\noindent since $g(R_i)=cR_i$. Hence $C_i=c(1-r)/(g'(R_i)-c)$ and,
using the formula for $\lambda(\xi)$, we obtain
$$
l(x)= (1-r)\delta(x)+\sum_{i=1}^nC_iR_ie^{-R_ix}
$$
and

\begin{eqnarray*}
r(u) & = & \int_u^\infty l(x)dx\\
& = & \sum_{i=1}^nC_ie^{-R_iu}\quad\textrm{for }u>0.
\end{eqnarray*}

\noindent This is en elegant generalization of
Cram\'{e}r-Lundberg's formula, which is obtained when only the
contribution from $R_1=R$ is included. Since $R_1<b_1<R_2<\ldots
<R_n$, we see that the first term $Ce^{-Ru}$ dominates, as
expected.

\subsubsection{An alternative derivation of Cram\'{e}r's formula for $r(u)$}

\noindent In the derivation of the formula for $r(u)$, the process
$U(t)$ was divided into successive upcrossings. This gave a
natural probabilistic interpretation of the quantities $r$ and
$k(y)$. The traditional method for determining $r(u)$ is to derive
an integral equation, that is well-known in renewal theory, and to
show that its solution is given by Cram\'{e}r's formula. Although
it does not provide the same insight concerning the probabilistic
structure of the solution, this method has the advantage of being
more direct. Also, it can be generalized to the case when $c$
depends on the value of $U(t)$, which is indeed a natural
extension. For the sake of completeness, we describe also this
analytic derivation.\medskip

\noindent We are looking for an equation for $r(u)$ as a function
of $u$, that is based on an analysis of what can happen in a small
interval $(0,h)$ just after $t=0$. Such equations are common in
the more general theory for Markov processes and are referred to
as \emph{backward equations}. There are basically two possible
scenarios that can occur in the interval $(0,h)$:

\begin{itemize}
\item[1.] With probability $e^{-\lambda h}\approx 1-\lambda h$ no
loss occurs. At time $h$ we then have $U(h)=-ch$ and ruin has not
yet occurred. Looking ahead from $h$, the ruin probability is
$r(u+ch)$, since the surplus has increased by $ch$ in the interval
$(0,h)$.
\item[2.] With probability $\approx \lambda h$ a loss occurs in
$(0,h)$. Let $x$ denote the amount of loss. If $x>u$, ruin
occurs immediately. If $x\leq u$, we have $U(h)\approx x$ and
ruin has not yet occurred. Looking ahead from $h$, the ruin
probability is $r(u-x)$, since the surplus has decreased  by $x$
in the interval $(0,h)$.
\item[(3.)] With probability $o(h^2)$, more than one loss occur in
$(0,h)$. As $h\rightarrow 0$, this possibility can be excluded.
\end{itemize}

\noindent Combining this gives
$$
r(u)=(1-\lambda h)r(u+ch)+\lambda h\int_0^ur(u-x)F(dx)+\lambda
h\int_u^\infty F(dx)+o(h^2)
$$
as $h\rightarrow 0$. If we assume that $r'(u)$ exists, we get the
equation
$$
cr'(u)-\lambda r(u)+\lambda\int_0^ur(u-x)F(dx)+\lambda
\bar{F}(u)=0,
$$
where $\bar{F}(u)=1-F(u)$. We will solve this equation with the
boundary condition that $r(u)\rightarrow 0$ as $u\rightarrow
\infty$ if $c>\lambda\mu$. A difficulty is that both $r(u)$ and
$r'(u)$ are included in the equation. However, it turns out that
$r(u)$ can be eliminated by partial integration of the third term.
Using the fact that $d\bar{F}(x)=-F(dx)$, we get
$$
\int_0^ur(u-x)F(dx)=r(u)-r(0)\bar{F}(u)-\int_0^ur'(u-x)\bar{F}(x)dx.
$$
Substituting this in the above equation yields
$$
cr'(u)=\lambda\int_0^u r'(u-x)
\bar{F}(x)dx+\lambda(r(0)-1)\bar{F}(u),
$$
Here $r(0)$ is a constant that can be determined from the boundary
condition. This is an equation that involves only $r'(u)$. To
solve it, introduce $r=\lambda\mu/c$ and
$k(u)=(1-F(u))/\mu=\bar{F}(u)/\mu$. The equation then becomes
$$
r'(u)=r\int_0^ur'(u-x)k(x)dx+r(r(0)-1)k(u).
$$
This is a renewal equation that, by a convolution operation, can
be written as
$$
r'(u)=r(k*r')(u)+r(r(0)-1)k(u).
$$
This equation can be solved by an iteration that converges when
$r<1$, that is, when $c>\lambda\mu$. We have

\begin{eqnarray*}
r'(u) & = & r(r(0)-1)\left(\sum_{m=0}^\infty
r^mk^{m*}\right)*k(u)\\
& = & (r(0)-1)\sum_{m=1}^\infty r^mk^{m*}(u).
\end{eqnarray*}

\noindent According to the boundary condition, $\int_0^\infty
r'(u)=r(\infty)-r(0)=-r(0),$ and, since $\int_0^\infty
k^{m*}(u)du=1$ for all $m$, we have

\begin{eqnarray*}
\int_0^\infty\sum_{m=1}^\infty r^mk^{m*}(u)du & = &
\sum_{m=1}^\infty r^m\\
& = & \frac{r}{1-r}.
\end{eqnarray*}

\noindent Using these relations, we obtain $-r(0)=(r(0)-1)\cdot
r/(1-r),$ that is, $r(0)=r$. Hence, since $l(u)=(1-r)\sum_0^\infty
r^mk^{m*}(u)$ and $k^{0*}(u)=\delta(u)=0$ when $u>0$, we get

\begin{eqnarray*}
r'(u) & = & -(1-r)\sum_{m=1}^\infty r^mk^{m*}(u)\\
& = & -l(u)\quad \textrm{for }u>0,
\end{eqnarray*}

\noindent By integration it follows that $r(u)=\int_u^\infty
l(y)dy$. This is the same formula that was derived in 3.2.3.

\subsubsection{Approximation of $r(u,t)$}

\noindent So far we have been concerned with
$r(u)=P(T(u)<\infty)$. However, it is also interesting to study
when ruin occurs if $T(u)<\infty$. In the following we will prove
a law of large numbers for $T(u)$, that states that, if
$T(u)<\infty$, then with large probability $T(u)\approx u\bar{t}$
as $t\rightarrow \infty$, where $\bar{t}$ is given by a formula
that includes the Lundberg exponent $R$. We will show that
exponential inequalities, analogous to the ones for $P(S(t)\geq
tx)$, hold for $P(T(u)\leq ut)$ when $t<\bar{t}$ and for $P(ut\leq
T(u)<\infty)$ when $t>\bar{t}$. The exponent can be expressed in
terms of the function $h(x)$.\medskip

\noindent In the derivation of Chernoff's inequality, $P(S(t)\geq
tx)\leq e^{-th(x)}$ for $x>\lambda\mu$, in Section 2.4.1, we
started with the relation E$\left[e^{\xi S(t)}\right]=e^{tg(\xi)}$
and picked a suitable value of $\xi$, depending on $x$. This
relation can be written as E$\left[e^{\xi S(t)-tg(\xi)}\right]=1$
for all $t$. We will first show that, for some $\xi$-values, this
relation holds also for the stochastic times $T(u)$ and $T(-u)$ so
that, for instance,

\begin{equation}\label{eq:wald}
\textrm{E}\left[e^{\xi S(T(u))-T(u)g(\xi)},T(u)<\infty\right]=1.
\end{equation}

\noindent Starting from this equation, which is called
\emph{Wald's identity}, we will derive inequalities for $T(u)$
analogous to the Chernoff bounds.\medskip

\noindent \emph{Proof of Wald's identity:}\smallskip

\noindent Since the process $S(t)$ has independent increments, for
$0<s<t$, we have that $S(t)-S(s)$ is independent of all events
$A_s$ and stochastic variables that concern the values of the
process up to time $s$. This implies that

\begin{eqnarray*}
\textrm{E}\left[e^{\xi S(t)-tg(\xi)},A_s\right] & = &
\textrm{E}\left[e^{\xi(S(t)-S(s))-(t-s)g(\xi)}\cdot e^{\xi
S(s)-sg(\xi)},A_s\right]\\
& = & \textrm{E}\left[e^{\xi(S(t)-S(s))-(t-s)g(\xi)}\right]\cdot
\textrm{E}\left[e^{\xi S(s)-sg(\xi)},A_s\right]\\
& = & \textrm{E}\left[e^{\xi S(s)-sg(\xi)},A_s\right],
\end{eqnarray*}

\noindent where the last equality comes from the fact that
$S(t)-S(s)$ has the same distribution as $S(t-s)$ and so its
generating function is $(t-s)g(\xi)$. Now let
$A_s=\{T(u)\in(s-ds,s]\}$. Note that when the values of the
process up to time $s$ are given, we can decide if $A_s$ has
occurred or not. We get

\begin{eqnarray*}
\textrm{E}\left[e^{\xi S(t)-tg(\xi)},T(u)\in ds\right] & = &
\textrm{E}\left[e^{\xi S(s)-tg(\xi)},T(u)\in ds\right]\\
& \approx & \textrm{E}\left[e^{\xi S(T(u))-T(u)g(\xi)},T(u)\in
ds\right].
\end{eqnarray*}

\noindent In the last equality we have used the facts that, since
$S(t)$ is right-continuous at the jump points, the difference
between $S(T(u))$ and $S(s)$ is at most $cds$ when $s-T(u)\leq
ds$, and that, with large probability, at most one jump occurs in
$ds$. Integrating the above relation for $s\in (0,t]$ yields
$$
\textrm{E}\left[e^{\xi S(t)-tg(\xi)},T(u)\leq t\right]=
\textrm{E}\left[e^{\xi S(T(u))-T(u)g(\xi)},T(u)\leq t\right].
$$
Recalling the definition of the Esscher-transformed distribution
of $S(t)$ in Section 2.4.2, we see that the left hand side can be
written as $P_\xi(T(u)\leq t)$, where $P_\xi(\cdot)$ is the
transformed measure. To establish Wald's identity we have to show
that this tends to 1 as $t\rightarrow \infty$. To this end,
remember that $S(t)$ is still a compound Poisson process under the
measure $P_\xi$, but the mean is changed to
E$_\xi[S(t)]=tg'(\xi)$. Hence E$_\xi[U(t)]=t(g'(\xi)-c)$. If $\xi$
is chosen so that this is strictly positive -- that is, so that
$g'(\xi)>c$ -- then, by the law of large numbers, $U(t)\rightarrow
\infty$ with $P_\xi$-probability 1 and it follows that
$P_\xi(T(u)<\infty)=1$. Since $g'(\xi)$ is an increasing function
of $\xi$, the condition that $g'(\xi)>c$ is fulfilled for
$\xi>\xi_c$, where $\xi_c$ satisfies $g'(\xi_c)=c$.\medskip

\noindent To summarize, we have showed that (\ref{eq:wald}) holds
for $\xi>\xi_c$, where $\xi_c$ is defined via the relation
$g'(\xi_c)=c$. Analogously it can be shown for $T(-u)$ that
$$
\textrm{E}\left[e^{\xi S(T(-u))-T(-u)g(\xi)},T(-u)<\infty\right]=1
$$
if $\xi<\xi_c$, since then the drift is strictly negative and
$P_\xi(T(-u)<\infty)=1$.\hfill$\Box$\medskip

\noindent From the picture of the definition of the Lundberg
exponent $R$ in Figure \ref{fig:s60} it can be seen that
$0<\xi_c<R$ if $c>\lambda\mu$ and $R<\xi_c<0$ if $c<\lambda\mu$.
We will now see how Wald's identity can be used to study $T(u)$
for $c>\lambda\mu$. First remember that $T(u)$ is defined as the
first time when $U(t)>u$, that is, the first time when
$S(t)>u+ct$. Together with Wald's identity this yields that, for
$\xi>\xi_c>0$,
$$
1\geq \textrm{E}\left[e^{\xi(u+cT(u))-
g(\xi)T(u)},T(u)<\infty\right],
$$
that is,
$$
e^{-\xi\mu}\geq \textrm{E}\left[e^{(c\xi-g(\xi))T(u)}
,T(u)<\infty\right].
$$
In particular, for $\xi=R$ we get
$$
P(T(u)<\infty)=r(u)\leq e^{-Ru},
$$
which is referred to as \emph{Lundberg's inequality}. This
inequality holds for all $u>0$ and the exponent is the same as in
the asymptotic approximation in Section 3.2.5.\medskip

\noindent The above inequality can be used to estimate $P(T(u)\leq
ut)$ (compare with the Chernoff bound from Section 2.4.1). If
$c\xi-g(\xi)\leq 0$, when $T(u)\leq ut$, we have
$(c\xi-g(\xi))T(u)\geq (c\xi-g(\xi))ut$ so that

\begin{eqnarray*}
e^{-u\xi} & \geq & \textrm{E}\left[e^{(c\xi-g(\xi))T(u)},T(u)\leq
ut\right]\\
& \geq & e^{(c\xi-g(\xi))ut}P(T(u)\leq ut),
\end{eqnarray*}

\noindent and hence

\begin{eqnarray*}
P(T(u)\leq ut) & \leq & e^{-u\xi-ut(c\xi-g(\xi))}\\
& = & e^{-ut((c+1/t)\xi-g(\xi))}.
\end{eqnarray*}

\noindent To get the best possible estimate we want to minimize
the exponent. This is done by picking $\xi$ such that
$g'(\xi)=c+1/t$ and, recalling the definition of $h(x)$ from
Section 2.4.1, we get
$$
P(T(u)\leq ut)\leq e^{-uth(c+1/t)}.
$$
The above calculations are valid under the assumption that
$c\xi-g(\xi)\leq 0$, that is, $R\leq \xi$, where $\xi$ is defined
by $g'(\xi)=c+1/t$. Since $g'(\xi)$ is strictly increasing, the
condition that $\xi\geq R$ is equivalent to $g'(\xi)\geq g'(R)$,
that is, to $1/t\geq g'(R)-c$. If we introduce $\bar{t}$, defined
by the relation $1/\bar{t}=g'(R)-c$, we see that the above
estimate holds for $t\leq \bar{t}$.\medskip

\noindent Analogously, if we pick $\xi$ such that $\xi_c<\xi$ and
$c\xi-g(\xi)\geq 0$ we can estimate $P(ut<T(u)<\infty)$. We get
$$
P(ut<T(u)<\infty)\leq e^{-ut((c+1/t)\xi-g(\xi))}.
$$
If $g'(\xi)=c+1/t$ and $\xi_c<\xi\leq R$, that is, if
$c=g'(\xi_c)<c+1/t\leq g'(R)$, or, equivalently,
$0<1/t<1/\bar{t}$, that is, $t\geq \bar{t}$, then it follows that
$$
P(ut<T(u)<\infty)\leq e^{-uth(c+1/t)}.
$$
For $\xi=R$ we get $t=\bar{t}$ and the exponent then becomes
$\bar{t}h(c+1/t)=\bar{t}((c+1/\bar{t})R-g(R))=R$, since
$cR=g(R)$.\medskip

\noindent To summarize, we have shown that there is a time
$\bar{t}$, defined by the relation $1/\bar{t}=g'(R)-c$, such that
the deviations from $u\bar{t}$ can be estimated by
$$
P(T(u)\leq ut)\leq e^{-uth(c+1/t)}\quad\textrm{for }t\leq \bar{t},
$$
and
$$
P(ut\leq T(u)< \infty)\leq e^{-uth(c+1/t)}\quad\textrm{for
}t\geq \bar{t},
$$
and $\bar{t}h(c+1/\bar{t})=R$ for $t=\bar{t}$.\medskip

\noindent We will soon see that the exponent $H(t):=th(c+1/t)$ is
a strictly convex function of $t$ with $\min_tH(t)=H(\bar{t})=R$.
This fact makes it possible to study $T(u)$ when $T(u)<\infty$.
Assume for example that $t<\bar{t}$. By Cram\'{e}r's
approximation, as $u\rightarrow \infty$ we then have

\begin{eqnarray*}
P(T(u)\leq ut|T(u)<\infty) & = & \frac{P(T(u)\leq ut)}{r(u)}\\
& \leq & \frac{e^{-uH(t)}}{r(u)}\\
& \approx & \frac{e^{-u(H(t)-R)}}{C}.
\end{eqnarray*}

\noindent This tends to 0 exponentially fast, since $H(t)>R$ for
$t<\bar{t}$. Analogously it can be seen that $P(ut\leq T(u)<
\infty|T(u)<\infty)$ tends to 0 exponentially fast when
$t>\bar{t}$ and $u\rightarrow\infty$. This has the following
important interpretation: When $t>\bar{t}$, the ruin probability
$r(u,ut)$ can be approximated by $r(u)\approx Ce^{-Ru}$, and when
$t<\bar{t}$, we have $r(u,ut)<<r(u)$, since $r(u,ut)\leq
e^{-uH(t)}$ and $r(u)\approx Ce^{-Ru}$ with $H(t)>R$.\medskip

\noindent \emph{Proof of the convexity of $H(t)$:}\smallskip

\noindent We have $H(t)=t((c+1/t)\xi-g(\xi))=t(c\xi-g(\xi))+\xi$,
with $g'(\xi)=c+1/t$. When $t$ varies, $\xi$ varies as well, and
we get

\begin{eqnarray*}
dH & = & (c\xi-g(\xi))dt+(tc-tg'(\xi)+1)d\xi\\
& = & (c\xi-g(\xi))dt.
\end{eqnarray*}

\noindent Hence $H'(t)=\frac{dH}{dt}=c\xi-g(\xi)$ and

\begin{eqnarray*}
H''(t) & = & (c-g'(\xi))\frac{d\xi}{dt}\\
& = & -\frac{1}{t}\cdot\frac{d\xi}{dt}
\end{eqnarray*}

\noindent From the equation for $\xi$ it follows that
$-dt/t^2=g''(\xi)d\xi$, that is,
$$
\frac{d\xi}{dt}=-\frac{1}{t^2g''(\xi)}<0.
$$
Thus $H''(t)>0$ and we have showed that $H(t)$ is strictly convex.
The function $H(t)$ attains its smallest value when
$H'(t)=c\xi-g(\xi)=0$, that is, when $\xi=R$ and $t=\bar{t}$. As
described above, we then have
$H(\bar{t})=\bar{t}(cR-g(R))+R=R$.\hfill$\Box$

\subsubsection{Approximation of $r(-u,t)$}

\noindent We will now show that the above estimates of $r(u,t)$
also hold for $r(-u,t)=P(T(-u)\leq t)$ with small modifications
when $c<\lambda\mu$. As we have seen, in this case we have
$R<\xi_c<0$ and it follows from Wald's identity that
$$
\textrm{E}\left[e^{\xi S(T(-u))-T(-u)g(\xi)},T(-u)<\infty\right]=1
$$
for $\xi<\xi_c$. An interesting difference as compared to the
previous case is that at the time of ruin we now have
$S(T(-u))=-u+cT(-u)$. This means that
$$
\textrm{E}\left[e^{-u\xi+(c\xi-g(\xi))T(-u)},T(-u)<\infty\right]=1.
$$
This is an equation for the generating function of $T(-u)$: Put
$w=c\xi-g(\xi)$ for $\xi<\xi_c$. Since
$\frac{dw}{d\xi}=c-g'(\xi)>0$ so that $\xi:=R(w)$ is uniquely
determined, this is a 1-1 relation. Hence
$$
\textrm{E}\left[e^{wT(-u)},T(-u)<\infty\right]=e^{uR(w)}.
$$
For $w=0$ we have $R(0)=R<0$ which gives the exact relation

\begin{eqnarray*}
P(T(-u)<\infty) & = & e^{Ru}\\
& = & e^{-|R|u},
\end{eqnarray*}

\noindent where $P(T(-u)<\infty)=r(-u)$.\medskip

\noindent As before we can also estimate $P(T(-u)\leq ut)$. If
$c\xi-g(\xi)\leq 0$ we have

\begin{eqnarray*}
e^{u\xi} & \geq & \textrm{E}\left[e^{(c\xi-g(\xi))T(-u)},
T(-u)<ut\right]\\
& \geq & e^{(c\xi-g(\xi))ut}P(T(-u)\leq ut)
\end{eqnarray*}

\noindent so that

\begin{eqnarray*}
P(T(-u)\leq ut) & \leq & e^{-ut(c\xi-g(\xi))+u\xi}\\
& = & e^{-ut((c-1/t)\xi-g(\xi))}\\
& = & e^{-uth(c-1/t)}
\end{eqnarray*}

\noindent if $\xi$ is chosen such that $g'(\xi)=c-1/t$. This is
possible if $g(\xi)\geq c\xi$, that is, if $\xi\leq R<0$ so that
$g'(\xi)=c-1/t\leq g'(R)$. Putting $1/\bar{t}=c-g'(R)$ gives the
condition $1/\bar{t}\leq 1/t$, that is, $t\leq \bar{t}$.
Analogously we obtain
$$
P(ut\leq T(-u)< \infty)\leq e^{-uth(c-1/t)}\quad \textrm{for
}t\geq \bar{t}.
$$
To summarize, we have the formulas $1/\bar{t}=c-g'(R)$,
$H(t)=th(c-1/t)$ and $r(-u)=e^{-|R|u}$. Furthermore,

\begin{eqnarray*}
P(T(-u)\leq ut|T(-u)<\infty) & = & \frac{r(-u,ut)}{r(-u)}\\
& \leq & e^{-u(H(t)+R)}\quad \textrm{for }t\leq \bar{t},
\end{eqnarray*}

\noindent and

\begin{eqnarray*}
P(ut\leq T(-u)< \infty|T(-u)<\infty) & = & \frac{r(-u)-r(-u,ut)}{r(-u)}\\
& \leq & e^{-u(H(t)+R)}\quad \textrm{for }t\geq \bar{t}.
\end{eqnarray*}

\noindent Since $H(t)$ is strictly convex with $H(t)\geq
H(\bar{t})=-R>0$, we can hence localize $T(-u)$ well near
$u\bar{t}$.

\subsubsection{An interpretation of the modified distribution $P_R(S(t)\in dx)$}

\noindent The Esscher transformed distribution $P_R(S(t)\in dx)$
is defined by
$$
P_R(S(t)\in dx)=e^{Rx-tg(R)}F(t,dx)
$$
and we have seen that

\begin{eqnarray*}
\textrm{E}_R\left[e^{\xi S(t)}\right] & = &
\textrm{E}\left[e^{(\xi+R)S(t)-tg(R)}\right]\\
& = & e^{t(g(\xi+R)-g(R))}.
\end{eqnarray*}

\noindent Under this measure, we have E$_R[S(t)]=tg'(R)$ and
E$_R[U(t)]=t(g'(R)-c)>0$ when $c>\lambda\mu$. Hence, by the law of
large numbers, $P_R(T(u)<\infty)=1$. Also, by the same theorem,
since E$_R[U(u\bar{t})]=u\bar{t}(g'(R)-c)=u$, we should expect
that $T(u)\approx u\bar{t}$ under the measure $P_R$ as
$u\rightarrow\infty$. As we have just seen, given that
$T(u)<\infty$, we have that $T(u)\approx u\bar{t}$ as
$u\rightarrow\infty$. Hence, it seems as if the measure $P_R$ gives
an approximate description of the conditional distribution of the
process $S(t)$ given that $T(u)<\infty$ as $u\rightarrow\infty$.
It is not hard to see that this is true for fixed $t$: Since
$P(t<T(u)<\infty|T(u)<\infty)\rightarrow 1$ as
$u\rightarrow\infty$, we have the relation

\begin{eqnarray*}
P(S(t)\in dx|T(u)<\infty) & = & \frac{P(S(t)\in
dx,T(u)<\infty)}{P(T(u)<\infty)}\\
& \approx & \frac{P(S(t)\in dx,t<T(u)<\infty)}{r(u)}.
\end{eqnarray*}

\noindent Because of the Markov property, this probability equals

\begin{eqnarray*}
\frac{P(S(t)\in dx)P(t<T(u)<\infty|S(t)=x)}{r(u)} & = &
\frac{P(S(t)\in dx)r(u-x+ct)}{r(u)}\\
& = & \frac{F(t,dx)r(u-x+ct)}{r(u)},
\end{eqnarray*}

\noindent since, if $S(t)=x$ we have $U(t)=x-ct$ and so the
surplus at time $t$ is $u-x+ct$. As $u\rightarrow \infty$ with $t$
fixed, we have $r(u)\approx Ce^{-Ru}$, implying that
$$
\frac{r(u-x+ct)}{r(u)}\rightarrow e^{Rx-Rct}.
$$
Hence the conditional distribution of $S(t)$ converges to
$F(t,dx)e^{Rx-Rct}$ and, since $g(R)=Rc$,

\begin{eqnarray*}
F(t,dx)e^{Rx-Rct} & = & F(t,dx)e^{Rx-tg(R)}\\
& = & P_R(S(t)\in dx).
\end{eqnarray*}

\noindent A corresponding result holds for $T(-u)$ when
$c<\lambda\mu$.\medskip

\noindent Using the distribution $P_R$ a fairly intuitive proof of
the central limit theorem for the quantity
$(T(u)-u\bar{t})/\sqrt{u}$ can be formulated as follows: If we
invert the relation between $P$ and $P_R$ we see that
$$
P(T(u)\leq u\bar{t}+t\sqrt{u})=\textrm{E}_R\left[
e^{-RS(T(u))+T(u)g(R)},T(u)\leq u\bar{t}+t\sqrt{u}\right].
$$
Furthermore, $U(T(u))=S(T(u))-cT(u)=u+Z$, where $Z$ is the
overshoot over $u$ at the passage at $T(u)$. The overshoot $Z$ is
bounded when $u$ is large and approximately independent of $T(u)$.
Hence, because $g(R)=cR$ the exponent in this expression can be
written $-R(cT(u)+u+Z)+cRT(u)=-Ru-RZ$, so that
$$
P(T(u)\leq u\bar{t}+t\sqrt{u})\approx
e^{-Ru}\textrm{E}_R\left[e^{-RZ}\right]P_R(T(u)\leq
u\bar{t}+t\sqrt{u}).
$$
The last probability can be estimated using the fact that, under
the modified measure $P_R$, the process $U(t)=S(t)-ct$ has
positive drift E$_R[U(t)]=t(g'(R)-c)=t/\bar{t}$ and variance
Var$U(t)=\textrm{Var}(S(t))=tg''(R)=t\sigma^2$. We can now
estimate $T(u)$ as follows. The law of large numbers tells us that
$U(t)/t\rightarrow 1/\bar{t}$ when $t\rightarrow \infty$ and, sice
$T(u)\rightarrow \infty$ as $u\rightarrow\infty$, we have
$U(T(u))/T(u)\rightarrow1/\bar{t}$ as $u\rightarrow\infty$. But,
since $U(T(u))=u+Z$ with $Z$ bounded, this implies that
$u/T(u)\rightarrow 1/\bar{t}$, that is, $T(u)/u\rightarrow
\bar{t}$. We can now use the central limit theorem for $U(t)$,
which tells us that the quantity
$$
X:=\frac{U(t)-t/\bar{t}}{\sigma\sqrt{t}}
$$
has an approximatively $N(0,1)$ distribution when
$t\rightarrow\infty$. Using this for $t=T(u)$ we get

\begin{eqnarray*}
X & = & \frac{U(T(u))-T(u)/\bar{t}}{\sigma\sqrt{T(u)}}\\
& \approx &
\frac{\bar{t}U(T(u))-T(u)}{\sigma\bar{t}^{3/2}\sqrt{u}},
\end{eqnarray*}

\noindent and, since $U(T(u))=u+Z$, this can be written
$$
X=\frac{\bar{t}u-T(u)+\bar{t}Z}{\sigma\bar{t}^{3/2}\sqrt{u}}.
$$
Since $Z$ remains bounded, $Z/\sqrt{u}$ can be neglected when
$u\rightarrow\infty$ and we finally get the Gaussian approximation
$$
\frac{T(u)-u\bar{t}}{\sqrt{u}}\approx -\sigma\bar{t}^{3/2}X
$$
under the measure $P_R$ and hence
$$
P_R\left(\frac{T(u)-u\bar{t}}{\sqrt{u}}\leq
\sigma\bar{t}^{3/2}x\right)\approx \Phi(x).
$$
Finally we get the corresponding formula for the measure $P$,
$$
P\left(\frac{T(u)-u\bar{t}}{\sqrt{u}}\leq\sigma\bar{t}^{3/2}x\right)\approx
C_Re^{-Ru}\Phi(x)
$$
with $C_R=\lim_{u\rightarrow\infty}\textrm{E}_R[e^{-RZ}]$. The
value of the constant $C_R$ can be deduced from the
Cram\'{e}r-Lundberg approximation $r(u)=P(T(u)<\infty)\approx
Ce^{-Ru}$ (see Section 3.2.5). If we let $x\rightarrow\infty$ we
see that $C_R=C=(c-g'(0))/(g'(R)-c)$. The asymptotic variance of
$T(u)$ is hence $u\bar{t}^3\sigma^2=ug''(R)/(g'(R)-c)^3$.

\subsubsection{An interesting property of a composite system}

\noindent Let us collect the approximate formulas for $r(u)$ and
$T(u)$ as follows. The exponent $R$ and the time $\bar{t}$ are
determined by $c=g(R)/R$ and $\bar{t}=1/(g'(R)-c)$. The
approximate time of ruin is $\bar{T}=u\bar{t}=u/(g'(R)-c)$ and, if
we define $C=(c-g'(0))/(g'(R)-c)$, then $r(u)\approx r(u,t)\approx
Ce^{-Ru}$ for $t>\bar{T}$. This means that, if our planning
horizon is $\bar{T}$, then the probability of ruin, $r(u)$, is a
reasonable approximation for the finite time ruin probability
$r(u,t)$ if $t>\bar{T}$. If ruin happens it takes place for
$T(u)\approx \bar{T}$.\medskip

\noindent Let us now consider a system consisting of two
independent pieces so that $S(t)=S_1(t)+S_2(t)$ with $S_1(t)$ and
$S_2(t)$ independent, and hence $g(\xi)=g_1(\xi)+g_2(\xi)$. It is
interesting to compare the quantities of the pieces to those of
the total system. It they have the same $R$, we get

\begin{eqnarray*}
c & = & \frac{g(R)}{R}\\
& = & \frac{g_1(R)}{R}+\frac{g_2(R)}{R}\\
& = & c_1+c_2,
\end{eqnarray*}

\noindent and, if they have the same $\bar{T}$, we obtain

\begin{eqnarray*}
u & = & \bar{T}(g'(R)-c)\\
& = & \bar{T}(g_1'(R)-c_1+g'_2(R)-c_2)\\
& = & u_1+u_2.
\end{eqnarray*}

\noindent If we use these $c_i$ and $u_i$, we get

\begin{eqnarray*}
r(u) & \approx & Ce^{-Ru}\\
& = & Ce^{-Ru_1}e^{-Ru_2}\\
& \approx & C_1e^{-Ru_1}C_2e^{-Ru_2}\\
& \approx & r_1(u_1)r_2(u_2),
\end{eqnarray*}

\noindent since, from the fact that

\begin{eqnarray*}
C & = & \frac{c-g'_1(0)+c_2-g'_2(0)} {g'_1(R)-c_1+g'_2(R)-c_2},
\end{eqnarray*}

\noindent it follows that $C_1\leq C_1C_2/C\leq C_2$ if $C_1\leq
C_2$, that is, the constants are comparable.\medskip

\noindent There is hence a natural decomposition of $c$ and $u$
into $c_1+c_2$ and $u_1+u_2$, so that if we have a common
$\bar{T}$ and $R$, then $r(u)\approx r_1(u_1)r_2(u_2)$, which is
the probability that both systems are ruined. None of the systems
is so to speak unnecessarily safe compared to the other. This is
also an example of decentralized planning: In order to calculate
$r(u)$ the central actuary only has to give the values of $R$ and
$\bar{T}$ to the local actuaries who can then calculate $r_1(u)$
and $r_2(u)$ and return them, and $r(u)\approx
r_1(u)r_2(u)$.

\vfill\eject

\section{Summary of the formulas}

\noindent In this section we give a concise summary of the
formulas that have been derived in the notes.\bigskip

\noindent \textbf{The individual risk model}\smallskip

\noindent $X$ = total amount of loss\newline\hspace*{0.35cm} =
$\sum_ix_iM_i$, where $\{M_i\}$ are Bernoulli variables with
$P(M_i=1)=p_i=1-q_i$.\medskip

\noindent Moments:
E$\left[X\right]=\sum_ix_ip_i$\newline\hspace*{1.5cm}
Var$(X)=\sum_ix_i^2p_iq_i$\medskip

\noindent Generating function: $\textrm{E}\left[e^{\xi
X}\right]=\prod_i\left(q_i+p_ie^{\xi x_i}\right)$\medskip

\noindent Compound Poisson approximation:\smallskip

\noindent $X\approx S$\newline \noindent $S=\sum_ix_iN_i$, where
$\{N_i\}$ are Poisson variables with $e^{-\lambda_i}=q_i$\medskip

\noindent Generating function: E$\left[e^{\xi
S}\right]=\exp\left\{\sum_i\lambda_i\left(e^{\xi
x_i}-1\right)\right\}$\newline\hspace*{4.4cm} $=e^{g(\xi)}$, where
$g(\xi)=\sum_i\lambda_i\left(e^{\xi x_i}-1\right)$\bigskip

\noindent \textbf{The collective risk model}\smallskip

\noindent $S(t)$ = total amount of loss in $(0,t)$\newline
\noindent $N(t)$ = number of accidents in $(0,t)$\newline\noindent
$X_i$ = the losses in the accidents\medskip

\noindent $S(t)=\sum_1^{N(t)}X_i$\medskip

\noindent $\{N(t)\}$ is a Poisson process with E$[N(t)]=\lambda
t$.\newline\noindent $\{X_i\}$ are i.i.d.\ with distribution
$F(dx)$, E$[X_i]=\mu$, E$[X_i^2]=\nu$.\newline\noindent The
distribution of $S(t)$ is $F(t,dx)=P(S(t)\in dx)$.\medskip

\noindent Generating function: E$\left[e^{\xi
S(t)}\right]=e^{tg(\xi)}$ with $g(\xi)=\lambda\int_0^\infty
\left(e^{\xi x}-1\right)F(dx)$.\medskip

\noindent Moments:
E$[S(t)]=tg'(0)=t\lambda\mu$\newline\hspace*{1.55cm}
Var$(S(t))=tg''(0)=t\lambda\nu$.\bigskip

\noindent \textbf{Panjer-recursion for the density of
$S(t)$}\smallskip

\noindent Assume that $X_i$ have a discrete distribution with
$P(X_i=nd)=f_n$. Then $P(S(t)=md)=g_m$ are given by the recursion
$$
\left\{\begin{array}{lll}
                     mg_m & = & \lambda t\sum_1^mnf_ng_{m-n},\quad m=1,2,\ldots \\
                     g_0 & = & e^{-\lambda t}.
                    \end{array}
            \right.
$$

\noindent \textbf{Approximations of $P(S(t)>tx)$}\smallskip

\noindent Entropy function: $
h(x)=\max_\xi\{x\xi-g(\xi)\}$\newline\hspace*{3.65cm}=
$x\xi_x-g(\xi_x)$, with $\xi_x$ defined by $g'(\xi_x)=x$.\medskip

\noindent The functions $g(\xi)$ and $h(x)$ are convex.\medskip

\noindent We have $g(\xi)=\max_x\{\xi
x-h(x)\}$\newline\hspace*{2.15cm}= $\xi x_\xi-h(x_\xi)$, with
$x_\xi$ defined by $h(x_\xi)=\xi$.\medskip

\noindent The functions $x=g'(\xi)$ and $\xi=h'(x)$ are inverses
of each other.\medskip

\noindent Chernoff's bound:\nopagebreak\smallskip
$$
\left\{\begin{array}{ll}
                    P(S(t)\geq tx)\leq e^{-th(x)} & \textrm{if } x\geq \lambda\mu;\\
                    P(S(t)\leq tx)\leq e^{-th(x)} & \textrm{if } x\leq \lambda\mu.
                    \end{array}
            \right.
$$

\noindent Esscher's approximation:\medskip

\noindent The Esscher transform of $F(dx)$ is
$F_a(dx)=e^{ax}F(dx)/f(a)$ with $f(a)=\int_0^\infty
e^{ax}F(dx)$.\medskip

\noindent E$_a\left[e^{\xi X_i}\right]=f(\xi+a)/f(a)$\newline
\noindent E$_a\left[e^{\xi S(t)}\right]=e^{t(g(\xi+a)-g(a))}=
e^{tg_a(\xi)}$\medskip

\noindent The transform of $F(t,dx)$ is $P_a(S(t)\in
dx)=F_a(t,dx)=e^{ax-tg(a)}F(t,dx)$.\medskip

\noindent Moments: E$_a[S(t)]=tg'(a)$\newline\hspace*{1.55cm}
Var$_a(S(t))=tg''(a)$\medskip

\noindent Esscher's approximation tells us that
$$
P(S(t)\geq tx)\approx \frac{e^{-th(x)}}{\sqrt{2\pi}a
\sqrt{tg''(a)}}
$$
with $x=g'(a)\geq \lambda\mu=g'(0)$, $a\geq 0$. This is valid for
a continuous distribution. For a discrete distribution with span
$d$, the factor $a$ is changed into $A(d)=(1-e^{-ad})/d$.\bigskip

\noindent \textbf{Ruin probabilities}\smallskip

\noindent $U(t)=S(t)-ct$ = net amount of loss in $(0,t)$\newline
\noindent $u$ = initial capital\newline \noindent
$T(u)=\min\{t;\hspace{0.1cm}U(t)>u\}$ = time of ruin\newline
\noindent $T(-u)=\min\{t;\hspace{0.1cm}U(t)=-u\}$\newline
\noindent $r(\pm u,t)=P(T(\pm u)\leq t)$ = ruin probabilities in
finite time\newline \noindent $r(\pm u)=P(T(\pm u)<\infty)$ = ruin
probabilities in infinite time\medskip

\noindent For $u=0$, we have the explicit formula
$$
P(T(0)>t,S(t)\in dx)=\left(1-\frac{x}{ct}\right)_+F(t,dx)
$$
and hence

\begin{eqnarray*}
P(T(0)>t) & = & 1-r(0,t)\\
& = & \int_0^{ct}\left(1-\frac{x}{ct}\right)F(t,dx).
\end{eqnarray*}

\noindent \textbf{The distribution of $T(-u)$}\smallskip

\noindent For $ct\geq u>0$ we have $P(T(-u)\in
dt)=(u/ct)F(t,cdt-u)$. Hence
$$
r(-u,t)=\int_{u/c}^\infty \left(\frac{u}{cs}\right)F(s,cds-u)
$$
and

\begin{eqnarray*}
r(-u) & = & \int_{u/c}^t
\left(\frac{u}{cs}\right)F(s,cds-u)\\
& = & \int_0^\infty
\left(\frac{u}{x+u}\right)F\left(\frac{x+u}{c},dx\right)\quad
\textrm{with }x=cs-u.
\end{eqnarray*}

\noindent If $c>\lambda\mu$, we have $r(-u)=1$.\bigskip

\noindent \textbf{The distribution of $T(u)$}\smallskip

\noindent Seal's formula:
$$
r(u,t)=\int_{u+ct}^\infty F(t,dx)+\int_0^tF(s,u+cds)\bar{r}(0,t-s)
$$
where $\bar{r}(0,t-s)=1-r(0,t-s)$.\bigskip

\noindent \textbf{Cram\'{e}r's formula for $r(u)$}\nopagebreak\smallskip

\noindent The upcrossings $\{U_k\}$ are i.i.d.\ with
$r:=P(U_1>0)=\lambda\mu/c$ if $c>\lambda\mu$, and $U_1$ has the
conditional density $k(u)=(1-F(u))/\mu$, given that
$U_1>0$.\medskip

\noindent The density of $\bar{U}=\max_{t\geq 0}U(t)$ is
$l(u)=(1-r)\sum_0^\infty r^mk^{m*}(u)$ and we have

\begin{eqnarray*}
r(u) & = & P(\bar{U}>u)\\
& = & \int_u^\infty l(y)dy.
\end{eqnarray*}

\noindent \textbf{Panjer-approximation of $r(u)$}\smallskip

\noindent Approximate the density $k(u)$ by a discrete one with
masses $\{k_n\}$ for $u=nd$, $n=1,2,\ldots$. Then the
corresponding approximation $\{l_n\}$ for $l(u)$, $u=nd$, can be
calculated by the iteration
$$
\left\{\begin{array}{lll}
                    l_n=r(k_1l_{n-1}+\ldots+k_nl_0)\quad \textrm{for }n\geq 1;\\
                    l_0=1-r.
                    \end{array}
            \right.
$$
The ruin probability $r(u)$ is approximated by $r_n=\sum_n^\infty
l_m$ for $u=nd$.\bigskip

\noindent \textbf{Cram\'{e}r-Lundberg's approximation of
$r(u)$}\nopagebreak\smallskip

\noindent For $c>\lambda\mu$, let $R$ be the positive root of the
equation $g(R)=cR$ and define $C=(c-g'(0))/(g'(R)-c)$. Then
$$
r(u)\approx Ce^{-Ru}\quad \textrm{as }u\rightarrow\infty.
$$
Similarly, for $c<\lambda\mu$, we have $r(-u)=e^{Ru}$, where $R$
is the negative root of $g(R)=cR$.\bigskip

\noindent \textbf{Approximation of $r(u,t)$}\smallskip

\noindent For $c>\lambda\mu$, define
$\bar{T}=u\bar{t}=u/(g'(R)-c)$. Then $T(u)\approx \bar{T}$ if
$T(u)<\infty$. More accurately, if $H(t)=th(c+1/t)$, then
$$
P(T(u)\leq ut|T(u)<\infty)\leq
C^{-1}e^{-u(H(t)-R)}\quad\textrm{for }t<\bar{t}
$$
and
$$
P(ut\leq T(u)<\infty|T(u)<\infty)\leq
C^{-1}e^{-u(H(t)-R)}\quad\textrm{for }t>\bar{t}
$$
and $R=\min_{t}H(t)=H(\bar{t})$. Similarly, for $c<\lambda\mu$, if
we define $\bar{T}=u\bar{t}=u/(c-g'(R))$ and $H(t)=th(c-1/t)$, we
have
$$
P(T(-u)\leq ut|T(-u)<\infty)\leq e^{-u(H(t)+R)}\quad\textrm{for
}t<\bar{t}
$$
and
$$
P(ut\leq T(-u)<\infty|T(-u)<\infty)\leq
e^{-u(H(t)+R)}\quad\textrm{for }t>\bar{t}
$$
and $-R=\min_tH(t)=H(\bar{t})$.\bigskip

\noindent \textbf{Interpretation of the transformed distribution
of $S(t)$}\smallskip

\noindent When $c>\lambda\mu$, the transformed distribution

\begin{eqnarray*}
F_R(t,dx) & = & P_R(S(t)\in dx)\\
& = & e^{Rx-tg(R)}F(t,dx)
\end{eqnarray*}

\noindent is equal to $\lim_{u\rightarrow\infty}P(S(t)\in
dx|T(u)<\infty)$ and the corresponding result holds for $T(-u)$
when $c<\lambda\mu$. Hence we have
$$
\textrm{E}[S(t)|T(u)<\infty]\rightarrow\textrm{E}_R[S(t)]=tg'(R)\quad\textrm{as
}u\rightarrow\infty.
$$
This explains the formula for $\bar{T}$, because
E$_R[U(t)]=t(g'(R)-c)$ so $\bar{T}$ is that value of $t$ for which
this is equal to $u$.\bigskip

\noindent \textbf{The central limit theorem for $T(u)$}\smallskip

\noindent When $u\rightarrow\infty$ we have
$$
P\left(\frac{T(u)-u\bar{t}}{\sqrt{u}}\leq\sigma\bar{t}^{3/2}x\right)
\approx Ce^{-Ru}\Phi(x),
$$
with $C=(c-g'(0))/(g'(R)-c)$, $\bar{t}=1/(g'(R)-c)$ and
$\sigma^2=g''(R)$.

\vfill\eject

\section{Notes and references}

\noindent The notes on risk theory by Harald Cram\'{e}r from 1930
[3] still form a very readable introduction to the subject. The
idea of using Lemma \ref{lemma:ruin} -- the so called \emph{ballot
theorem} -- to derive the formulas for ruin probabilities is
developed by Lajos Tak\'{a}cs in [6]. Hopefully our treatment is
more understandable. The use of tools from large deviation theory
to derive asymptotic estimates is developed by the author in [5].
An alternative way of studying $T(u)$, which allows a central
limit theorem to be proved is developed by Bengt von Bahr in [2].
A modern and comprehensive treatment of the theory of ruin
probabilities is given in [1] and [4].

\vspace{0.5cm}

\begin{itemize}

\item[$\mbox{[1]}$ ]Asmussen, S: \emph{Ruin probabilities}, World Scientific 2000.

\item[$\mbox{[2]}$ ]von Bahr, B: Ruin probabilities expressed in terms of ladder
height distributions, \emph{Scand. Actuarial J.} 1974, 190-204.

\item[$\mbox{[3]}$ ]Cram\'{e}r, H: On the mathematical theory of risk,
\emph{H.C. Collected Works} vol. \textbf{1}, 601-678, Springer
1994.
\item[$\mbox{[4]}$ ]Klugman, S, Panjer H, Willmot, G: \emph{Loss Models, 2 ed.}, Wiley 2004.

\item[$\mbox{[5]}$ ]Martin-L\"{o}f, A: Entropy a useful concept in risk
theory, \emph{Scand. Actuarial J.} 1986, 223-235.

\item[$\mbox{[6]}$ ]Tak\'{a}cs, L: \emph{Combinatorial methods in the
theory of stochastic processes}, Wiley 1967.

\end{itemize}

\end{document}